\documentclass[12pt,twoside]{article}
\usepackage{graphicx}
\usepackage{amsfonts}
\usepackage{amsbsy}
\usepackage{amssymb}
\usepackage{amsmath, amsthm}
\usepackage{cite}
\usepackage{pst-all}
\usepackage{pstricks}
\usepackage{graphics}
\usepackage{tikz-cd}
\usepackage{bm, enumitem}

\def\bpsp{\begin{pspicture}}
\def\epsp{\end{pspicture}}

\newcommand{\R}{\mathbb{R}}
\newtheorem{theorem}{Theorem}[section]

\newtheorem{example}[theorem]{Example}
\newtheorem{lemma}[theorem]{Lemma}

\newtheorem{definition}[theorem]{Definition}
\newtheorem{proposition}[theorem]{Proposition}
\newtheorem{note}{Note}

\newtheorem{prop}[theorem]{Proposition}

\newcommand{\fos}{\mbox{ for some }}

\newtheorem{alg}[theorem]{Algorithm}
\newtheorem{notations}[theorem]{Notations}
\newtheorem{examples}[theorem]{Examples}
\newtheorem{conjecture}[theorem]{Conjecture}
\newlist{myenu}{enumerate}{5}
\setlist[myenu,1]{label=\arabic*.}
\setlist[myenu,2]{label*=\arabic*.}
\setlist[myenu,3]{label*=\arabic*.}
\setlist[myenu,4]{label*=\arabic*.}
\setlist[myenu,5]{label*=\arabic*.}

\newcommand{\wimp}{,\mbox{ which implies }}
 
\newcommand{\Nt}{Note that }
\newcommand{\nt}{note that }
\newcommand{\st}{such that }

\newcommand{\wv}{we have }
\newcommand{\wo}{we obtain }
\newcommand{\aThr}{Therefore }

\newcommand{\Csq}{Consequently, }

 \newcommand{\sm}{\setminus}
 
 \newcommand{\T}{\mathcal{T}}
  \newcommand{\hst}{\hspace{.2in}}

 \newcommand{\tcal}{\mathcal{T}}

 \newcommand{\nd}{\mbox{ and }}

  \newcommand{\foa}{\mbox{ for all }}
\newcommand{\ndfoa}{\mbox{ and for all }}
 \newcommand{\del}{\delta}
 \newcommand{\ds}{\displaystyle}

\begin{document}
\title{Sphere of Influence Dimension Conjecture 'Almost Proved'}
\author{Surinder Pal Singh Kainth$^a$,  Ramanjit Kumar$^b$, S. Pirzada$^c$\\
$^{a}${\em Department of Mathematics, Panjab University, Chandigarh, India}\\
$^{b}${\em 6, Ganga Enclave, Kot Khalsa, Amritsar, Punjab, India}\\
$^{c}${\em Department of Mathematics, Kashmir University, India}\\
   $^a$sps@pu.ac.in;  $^b$ramanjit.gt@gmail.com; $^c$pirzadasd@kashmiruniversity.ac.in}
\date{}
\pagestyle{myheadings}
\markboth{S. P. S. Kainth, R. Kumar and S. Pirzada}{SIG-Dimension Conjecture}
\maketitle
\vskip 5mm
\noindent{\footnotesize \bf Abstract.}
The sphere-of-influence graph (SIG) on a finite set of points in  a metric space, each with an open ball centred about it of radius equal to the distance between that point and its nearest neighbor, is defined to be the intersection graph of these balls. Let $G$ be a graph of order  $n$ having no isolated vertices. The SIG-dimension of  $G,$ denoted by $SIG(G),$ is defined to be the least possible $d$ such that $G$ can be realized as a sphere of influence graph in $\R^d,$ equipped with sup-norm. In 2000, Boyer [E. Boyer, L. Lister and B. Shader, Sphere of influence graphs using the sup-norm, Mathematical and Computer Modelling 32 (2000) 1071-1082] put forward the SIG dimension conjecture, which states that
$$SIG(G)\leq  \bigg\lceil \frac{2n}{3}\bigg\rceil.$$ 
In this paper, we 'almost' establish this conjecture by proving that 
$$SIG(G)\leq \bigg{ \lfloor}\frac{2n}{3}\bigg{ \rfloor}+2.$$
%

\vskip 3mm

\noindent{\footnotesize Keywords: Sphere of influence, matching, factor,  star-triangle factor.}

\vskip 3mm
\noindent {\footnotesize AMS subject classification: Primary: 05C62, 05C75, 05C70; Secondary:  68R10.}



\section{\bf Introduction}

A simple \textit{graph} is denoted by $G(V(G),E(G))$, where $V(G)=\{v_{1},v_{2},\ldots,v_{n}\}$ is its vertex set and $E(G)$ is its edge set. The \textit{order} of $G$ is $|V(G)|$ and its \textit{size} is $|E(G)|$.
A cycle of order $n$ is denoted by $C_{n}$ and a \textit{triangle} is $C_{3}$. A complete bipartite graph $K_{1,n}$ is called a \textit{star.} In $K_{1,n}$, the vertex of degree $n$ is its \textit{center} and all other vertices are \textit{leaves.}

An \textit{isolated vertex} of a graph is a vertex which is of degree zero.
Two edges  in  a graph are said to be \textit{independent edges} if  they do not share a common vertex. A \textit{matching} in a graph is a set of independent edges. That is, a subset $M$ of the edge set $E$ of $G$ is a matching if no two edges of $M$ have a common vertex. A \textit{perfect matching} is a matching which covers all vertices of the graph.

A matching $M$ is said to be \textit{maximal} if there is no matching $N$ strictly containing $M$, that is, $M$ is maximal if it cannot be enlarged. A matching $M$ is said to be \textit{maximum} if it has the largest possible cardinality, that is, $M$ is maximum if there is no matching $N$ such that $|N|~>~|M|$. A vertex $v$ is said to be \textit{$M$-saturated} (or saturated by $M$) if there is an edge $e\in M$ incident with $v$. A vertex which is not incident with any edge of $M$ is said to be \textit{$M$-unsaturated.} An \textit{$M$-alternating path} in $G$ is a path whose edges are alternately in $E(G)\sm M\nd M$. That is, in an $M$-alternating path, the edges alternate between $M$-edges and non-$M$-edges. An $M$-alternating path whose end vertices are $M$-unsaturated is said to be an \textit{$M$-augmenting path.}

We construct the   sphere of influence graphs with sup-norm as follows.
 Let $d$ be a natural number and $\ensuremath{\mathbb{R}}^d$ denote the $d$-dimensional Euclidean space. For any $z\in \ensuremath{\mathbb{R}}^d,$ let $z[j]$ denote the $j^{th}$ component of $z$. The distance between any $x, y\in \ensuremath{\mathbb{R}}^d$ under the $L_\infty$-metric, denoted by $\rho(x,y),$ is defined as
  $$\rho(x,y):=\max \{|x[j]-y[j]|: j=1,2,\dots,d\}.$$
Let $P\subset \ensuremath{\mathbb{R}}^d $ be a finite  set having at least two points. For a point $v\in P,$ let
$r_v$ denote the distance of $v$ to its \textit{nearest neighbor,} that is
$$r_v=\min\{\rho(u, v):u\in P\setminus \{v\}\}.$$
The open ball $B_v:=\{u\in \ensuremath{\mathbb{R}}^d:\rho(u,v)<r_v\}$ is known as the \textit{sphere of influence at $v.$}
The \textit{sphere of influence graph of $P,$} denoted by $SIG_\infty^d(P),$ is the graph with vertex set $P$ and edges corresponding to the pairs of intersecting spheres of influence. That is, the edge set of $SIG_\infty^d(P)$ is given by
$$\{uv: B_u\cap B_v\neq \emptyset; u, v \in P\}.$$

For $G=SIG_\infty^d(P)$ and  $u, v \in P,$ we observe that $uv\in E(G)$ if and only if $\rho(u,v)<r_{u}+r_{v}$.
A graph $G$ is said to be \textit{realizable} in $\ensuremath{\mathbb{R}}^d$ if there exists a finite set $P\subset \ensuremath{\mathbb{R}}^d$ such that $G$ is isomorphic to
$SIG_\infty^d(P).$ Note that if $G$ is {realizable} in $\ensuremath{\mathbb{R}}^d,$ then it is realizable in $\ensuremath{\mathbb{R}}^{d+e}$ for every $e\in \ensuremath{\mathbb{N}}.$ This can be observed by appending $e$ zero-coordinates to each point in the vertex set. The smallest such $d$ is called the SIG-dimension of the graph $G,$ denoted by $SIG(G).$ That is,
$$SIG(G)  = \min \{d : G\mbox{ is realizable  in }\ensuremath{\mathbb{R}}^d\}.$$

Note that if a graph with at least two vertices is realizable in some $\ensuremath{\mathbb{R}}^d,$ then it is trivial to see that it has no isolated  vertices.  Also,  every graph $G$ with at least two vertices and without isolated  vertices can be realized  in  $\ensuremath{\mathbb{R}}^d,$ for some $d\in \ensuremath{\mathbb{N}},$ as the rows of the matrix  $2I+  A$ realize $G$, where $A$ is the adjacency  matrix  for $G\nd I$ is the identity  matrix. For more details see \cite[Theorem 1]{2003}.

The notion of the sphere of influence graphs was introduced by Toussaint to model situations in pattern recognition and computer vision. These are used to help separate objects or otherwise capture perceptual relevance, see \cite{1980, 1982, 1998}. Toussaint has used the sphere-of-influence graph  (SIGs) under $L_2$-norm  to capture low-level perceptual information in certain dot patterns. The SIGs in general metric spaces are discussed in \cite{1999}. It is known that the SIGs under the $L_\infty$-norm perform better for this purpose, see \cite{2003}.

The \textit{SIG-dimension} of a graph $G$, denoted by $SIG(G)$,  is the minimum possible $d$ such that $G$ can be realized as a  sphere of influence graph in $\R^d$ with $L_\infty$-norm. Recently in \cite{2014}, Taussaint has surveyed the theory and applications of sphere of influence graphs. In \cite{2003}, several open problems on SIG-dimension have been discussed and the one regarding SIG-dimension of trees has already been solved, for details see \cite{2011}. Boyer et al. \cite{2000} put forward the following conjecture.
\begin{conjecture}[Boyer et al., 2000]\label{conj}
If a graph $G$ of order $n$ has no isolated vertices, then
$$SIG(G)  \leq \bigg{ \lceil}\frac{2n}{3}\bigg{ \rceil}.$$
\end{conjecture}
For more on this conjecture we refer Section 5 of \cite{2000}, where   it has been verified  for a few  families of graphs. In \cite{rs2}, Kumar and Kainth have established  this conjecture for graphs having  perfect matchings. This also generalizes various available results in this direction, see corollaries 17 and 18 of \cite{2000}. In \cite{rs2}, the following weaker bound for the SIG-dimension has also been established, for general graphs.
\begin{theorem}
If $G$ is a graph of order $n$ without isolated vertices, then
$$SIG(G)  \leq \bigg{ \lceil}\frac{3n}{4}\bigg{ \rceil}+\lceil \log_2 n\rceil.$$
\end{theorem}
Kumar and Kainth \cite{rs1} have also obtained a bound for the SIG-dimension in case of $K_{2,2}$-free graphs.
Recall that a graph is said to be $K_{2,2}$-free, if it does not contain an induced subgraph isomorphic to $K_{2,2}.$
They proved that if $G$ is a  $K_{2,2}$-free graph of order $n,$ without an isolated vertex, then
$SIG(G)  \leq \big{ \lfloor}\frac{3n}{4}\big{ \rfloor}+\lceil \log_2 n \rceil +2$.
A few partial results regarding the SIG-dimension for some particular graphs are proved in \cite{2000,  2003}.
In this paper, we 'almost' prove this conjecture, by  establishing  the following.

\begin{theorem}\label{MainThmSIG3}
 If $G$ is any graph of order $n$ and has no  isolated vertex, then
$$SIG(G)\leq \bigg{ \lfloor}\frac{2n}{3}\bigg{ \rfloor}+2.$$
In particular, if $n$ is not a multiple of $3,$ then
\begin{equation*}\label{MainThmMainEq}
SIG(G)\leq \bigg{ \lceil}\frac{2n}{3}\bigg{ \rceil}+1.
\end{equation*}
\end{theorem}

If $S\subset V(G),$ the \textit{induced graph} on $S$ is a  subgraph  of $G$ with vertex set $S$ and edge set consisting of all the edges of $G$ which have both end vertices in $S.$ An \textit{induced star} of $G$ is an induced subgraph of $G$ which itself is a star.

For a set $S$ of connected graphs, a spanning subgraph $F$ of a graph $G$ is called an \textit{$S$-factor} of $G$ if each component of $F$ is isomorphic to an element of $S$. A spanning subgraph $F$ of a graph $G$ is a \textit{star-cycle} factor of $G$ if each component of $F$ is a star or a cycle. A spanning subgraph $S$ of a graph $G$ will be called an \textit{induced star-triangle factor} of  $G$  if each component of $S$ is either an induced  star ($K_{1,n},~n\geq 2$, or $K_{2}$) or a triangle of $G$.


In \cite{krp},  the authors have established that every graph without isolated vertices admits an induced star-triangle factor in which   any two leaves from different stars $K_{1,n}$ ($n\geq 2$) are non adjacent.
This induced star-triangle factor will be used  here to embed a given graph into a suitable Euclidean space, in order to prove our main results. Further, through out we consider the graph $G$ with $n$ vertices and having no isolated vertices.

 \section{\bf The Main Ideas}

This paper is based  upon two main algorithms. The first one provides an induced star triangle factor of our given graph $G,$ in which leaves of no two stars are adjacent.  The second one provides a systematic way of picking up subsets of vertices of $G,$  which facilitates our embedding of graph $G$ into a suitable Euclidean space.

This section is divided into six further subsections. The first two are devoted to the two algorithms, mentioned above.
The third subsection provides some notations and definitions to be used again and again, throughout this manuscript.
In the fourth subsection, we shall classify our vertices of $G,$ based upon the second algorithm.
The fifth subsection will briefly outline the strategy to embed our graph into a suitable Euclidean space. In the last subsection, we shall discuss some sample inequalities, which will be used throughout the rest of the manuscript.

\subsection{Obtaining a Star-Triangle Factor}

Eager readers  may go through the next paragraph and then jump to subsection \ref{charversec}.
Others may just skip the next paragraph.

As established in \cite{krp}, $G$ admits an induced star-triangle factor in which  any two leaves from different stars $K_{1,n}$ ($n\geq 2$) are non adjacent.  If  $ K_{1,n},$ where $ n\geq 2$ is such a star  with center $u,$ we will write $S_u$ for the set of its leaves, that is, $K_{1,n}\sm \{u\}.$ Let $\mathcal{T}$ be the collection of all the triangles and $M$ be the collection of stars of the form $K_{1,1} $ of this induced star-triangle factor of $G.$ Note that $M$ is  actually a matching in $G.$

To have an insight into the main results of \cite{krp}, let $M_0$ be the maximum matching in $G.$ Let $M_0'$ be the set of matched vertices ($M$-saturated) and $I$ be the set of unmatched vertices ($M$-unsaturated).
Then $I$ must be an  independent set.

As in \cite{krp}, the following algorithm  provides an induced star triangle factor of our given graph $G,$ in which leaves of no two stars are adjacent.

\begin{alg}\label{mainalg}
 \begin{enumerate}
\item Let $M = M_0.$
\item \label{a2} If $I \ne \phi,$ then pick a vertex $v$ from $I$, otherwise go to step \ref{a10}.
\item Pick $u\in V(G)$  such that  $uv \in E(G)$ and call this edge the smallest edge of $v.$
(As $v$ is not isolated, there exists an edge $uv \in E(G).$)
Then $u\in M'_0.$ Otherwise, $M_0\cup \{uv\}$  will be a larger matching than $M_0.$
\item Choose $w\in M'_0$ such that $uw\in M_0.$
\item If $S_u$ is not defined, define $S_u:= \{w, v\}$, otherwise go to step \ref{a6}.
\item Remove $uw$ from $M,$ go to step \ref{a8}.
\item \label{a6} If $S_u$ is defined, then add $v$ to $S_u.$
\item \label{a8} Set $J = I \setminus \{v\}.$
\item With $I = J$, go to step \ref{a2}.
\item \label{a10} Stop.
\end{enumerate}\end{alg}
At the end of algorithm, we have $S_{u_1}, S_{u_2}, \ldots, S_{u_k}\nd u_1, u_2, \ldots, u_k\nd M.$
 The following results have been established in \cite{krp}.

\begin{prop} The residue set $M$ is a matching
and if  $M'$ is the set of vertices of $M,$ then $V(G)$ can be partitioned as
\begin{equation*}
V(G)=\big(\dot\cup_{S_u}(\{u\}  \cup S_u)\big)\dot \cup  M^{'} .
\end{equation*}
\end{prop}

\begin{proposition}
 If  $S_u$ is defined for some $u\in V(G)$ and for some $v,w\in S_u,$ $vw\in E(G),$ then $$S_u=\{w,v\}.$$
\end{proposition}

\noindent  We now make a little change in our notations from Algorithm \ref{mainalg}.
\begin{note}\label{note_1}
For each $S_u = \{v_1, v_2 \}$,  if $v_1v_2 \in E(G)$, then destroy (remove) $S_u.$
We mean, from now onwards this $S_u$ does not exist. Instead, if such an $S_u$ exists,
we do the following:

If $\T$ is not defined,   define $\T:= \{\{u, v_1, v_2\}\}$,
otherwise add $\{u, v_1, v_2\}$ to $\T.$
\end{note}

\begin{prop}\label{p3.8}
The set $\cup S_u$ is  independent.
\end{prop}

Summarizing all these observations,   the following result is established in \cite{krp}.

\begin{theorem}\label{mainthmsig3a}
Every  graph without isolated vertices, admits an induced star-triangle factor in which  any two leaves from different stars $K_{1,n}$ ($n\geq 2$) are non adjacent.
\end{theorem}


\subsection{Characterizing Vertices}\label{charversec}

To establish our main result, we need to realize our given graph $G$ in a suitable Euclidean space.
As a first step,
 we will be picking up vertices of the given graph $G,$ inductively, based upon the procedure specified in the following algorithm. It must be noted that the vertices once picked, won't be picked up again.

  The following  is the second major algorithm to be used in this paper.





\begin{alg}\label{alg41}
\begin{enumerate}
\item Let $A=\{S_u:S_u$ exists$\}\nd B=\emptyset.$
\item\label{CharVerS02} If $|A|+|B|<3,$ go to step \ref{CharVerS18}.
\item If $|B|=0$,  pick $S_x,S_y\nd S_z$ from $A.$
\item If $|B|=1$,  pick $S_x,S_y$ from $A\nd S_z$ from $B.$
\item If $|B|=2$,  pick $S_x$ from $A\nd S_y,S_z$ from $B.$
\item If $|B|=3$,  pick $S_x,S_y\nd S_z$ from $B.$
 \item\label{CharVerS07} If $\{x,y, z\}$ is an independent set, then pick $x,y,z$ and go to step \ref{CharVerS11}.
\item Let $x_1 \in S_x,y_1\in S_y\nd z_1 \in S_z.$
\item Out of the sets of triplets $\{x_1,y, z\}\mbox{ or }\{x,y_1, z\} \mbox{ and }\{x,y, z_1\},$ pick the  independent ones
and go to step \ref{CharVerS11}.
\item Out of the sets of triplets $\{x,y_1, z_1\}\mbox{ or }\{x_1,y, z_1\}\nd \{x_1,y_1, z\},$
pick the  independent ones
and go to step \ref{CharVerS11}.
\item Pick $x_1, y_1\nd z_1.$ (By  Proposition \ref{p3.8}, $x_1, y_1\nd z_1$ are independent.)
 \item\label{CharVerS11} If $v$ is picked, $S_v$ exists  and $S_v\in A,$ then remove $S_v$ from $A.$
\item If $v$ is picked, $S_v$ exists and $S_v\in B,$ then remove $S_v$ from $B.$
\item If $v$ is picked and $v\in S_{u},$ for some $u\ne v,$ then remove $v$ from $S_{u}.$
\item If $S_{u}\ne \phi\nd S_{u}\in A,$ then add it to $B.$
\item If $S_{u}=\phi\nd S_{u}\in B,$ then remove $S_{u}$ from $B.$
\item\label{CharVerS17} Go to step \ref{CharVerS02}.
\item\label{CharVerS18} If $|A|+|B|=2,$ then write $A\cup B=\{S_x,S_y\}$ and pick $x\nd y.$
\item\label{CharVerS19} If $|A|+|B|=1,$ then write $A\cup B=\{S_x\}$ and pick $x.$
\item Let $C=V(G)\setminus \{v:v$ is picked$\},$ be the set of unpicked vertices of $G.$
(Note that if $u\in V(G)$ is such that $S_u$ was defined and $u$ is picked, then $S_u \neq \emptyset.$)

\item Let $D_u=S_u$ if $u$ is picked.
\item \label{CharVerS22} If there are three vertices $p,q,s\in C$ such that there is no edge among them
and no two of $\{p,q,s\}$  are in same $D_u,$ then  pick $p,q,s.$ Otherwise,
go to step \ref{CharVerS24}.
\item Set $C=C\setminus \{p,q,s\}$ and go to step \ref{CharVerS22}.
\item \label{CharVerS24} \label{s30} If there are three vertices
$p,q,s\in C$ such that there is exactly one edge among them
and no two of $\{p,q,s\}$  are in same  $D_u,$ then  pick $p,q,s.$ Otherwise,
go to step \ref{CharVerS26}.
\item Set $C=C\setminus \{p,q,s\}$ and go to step \ref{s30}.
\item \label{CharVerS26} Note that we have at most two $D_u$ which have at least one unpicked vertex.
If there is no such $D_u,$ then go to step \ref{CharVerS40}.
\item\label{CharVerS27} Else, if there are exactly two such $D_u,$ then pick all unpicked vertices of these $D_{u}'s$ and call them $\{a_1,\dots, a_m\}.$
\item Set $C=C\sm \{a_1,\dots, a_m\}$ and go to step \ref{CharVerS40}.
\item Note that  there is  exactly one such $D_u,$ which has at least one unpicked vertex. Call it $D_w.$
\item \label{CharVerS30} If $|D_w|\leq 2,$ pick the whole of $D_w.$ Else go to step \ref{CharVerS32}.
\item Set $C=C\sm D_w$ and go to step \ref{CharVerS40}.
\item \label{CharVerS32} If there is no unpicked vertex, apart from the unpicked vertices of  $D_w,$ pick all these vertices of $D_w,$ call them $\{b_1,\dots,b_l\}$ and go to  step \ref{CharVerS47}. Else, go to step \ref{CharVerS33}.
\item \label{CharVerS33} If there is some unpicked vertex $v_0$ outside  $D_w\nd w_1,w_2\in D_w$ \st $\{v_0,w_1,w_2\}$ forms an independent set, pick this triplet. Else, go to step \ref{CharVerS35}.
\item Set $C=C\sm \{v, w_1,w_2\}\nd D_w=D_w\sm \{w_1,w_2\}.$
 Go to step \ref{CharVerS30}.
\item \label{CharVerS35}
If there are $v_1, v_2\in C\sm D_w$ with $v_1v_2\in E(G);\nd w_0\in D_w$ \st $w_0$ has an edge with exactly one of $ v_1\nd  v_2,$  pick $\{w_0, v_1, v_2\}.$ Else go to step \ref{CharVerS38}.
\item Set $C=C\sm \{w_0, v_1, v_2\}\nd D_w=D_w\sm \{w_0\}.$
Go to step \ref{CharVerS30}.
\item \label{CharVerS37}
If $C\sm D_w$ is a singleton set, pick that singleton set.
\item \label{CharVerS38}
Pick the whole of $D_w,$ say $\{c_1,\dots,c_{l'}\}.$
\item Set $C=C\sm D_w.$
\item \label{CharVerS40}If there are two vertices
$p,q\in C$ such that $pq\notin E(G)$ then pick $p,q.$ Else, go to step \ref{CharVerS42}.
\item Set $C=C\setminus \{p,q\}$ and go to step \ref{CharVerS40}.
\item \label{CharVerS42}  Note that the set $C,$ of unpicked vertices, forms a clique (a complete graph).
\item \label{CharVerS43} If $|C|\geq 3,$ then pick three vertices
$p',q',s'\in C,$ otherwise go to step \ref{CharVerS45}.
\item Set $C=C\setminus \{p',q',s'\}$ and go to step \ref{CharVerS43}.
\item \label{CharVerS45} If $|C|=1$ or $2,$ then pick the whole of $C.$
\item \label{CharVerS47}  Stop.
\end{enumerate}
\end{alg}
Through this algorithm, we have picked up sets of vertices of the form
\begin{gather*}
A,
X, \{p,q,s\}, \{a_1,\dots, a_m\}, Y, \{b_1,\dots,b_l\},\{v_0,w_1,w_2\},\\
\{w_0, v_1, v_2\},  \{c_1,\dots,c_{l'}\}, \{p, q\}, \{p',q',s'\} \nd  Z,
\end{gather*}
where
each of $X,Y \nd Z$  are  either  singleton sets or   pairs of vertices of $G\nd A$ is a triplet consisting of independent vertices, which are either the centers or the leaves of some stars.

Let us denote all these  picked up sets as, $P_0,P_1,\dots,P_{i},\dots,$ in the order they appear in the algorithm.
It is pertinent to mention that if for a particular graph a triplet of the form $\{x_{i},y_{i},z_{i}\}$ is picked three times, then we will denote those three with different notations $P_0, P_1\nd P_2.$

Note that we have used different letters of the alphabet to denote chosen vertices at different stages.
It is done only  for notational convenience. As we will be proceeding with more technical stuff, different notations will certainly help in recognizing the nature of the corresponding vertices.



\subsection{Notations and Definitions}

Before we move ahead, we fix some notations to be used in the sequel.
First we introduce the sets of pseudo-neighbors and pseudo-neighbor of pseudo-neighbors of the vertices of $G.$
\begin{definition}[Pseudo-neighbors]\label{nndefn}
 \begin{myenu}
\item Corresponding to each   $S_u,$  we define $$N(u):=S_u\mbox{ and  }N(x):=\{u\},\mbox{ for each }x\in S_u.$$
\item For each triangle $\{x,y,z\}\in \T,$ let $x\in P_{m_1}, y\in P_{m_2}\nd z\in P_{m_3}.$
Without loss of generality, let $m_1\le m_2\le m_3.$
We define $$N(x):=\{y,z\}, N(y):=\{x\}\mbox{ and }N(z):=\{x\}.$$
\item For $uv\in M,$ we define $N(v):=\{u\}$ and  $N(u):=\{v\}.$
\end{myenu}\end{definition}
Next, we define pseudo-neighbors of pseudo-neighbors of vertices of $G,$ naturally.
\begin{definition}\begin{enumerate}
  \item For any $v\in V(G),$ we define
 $$N^2(v):=\bigcup_{u\in N(v)}N(u)\nd N^2_0(v):= N^2(v)\sm \{v\}.$$
 \item A vertex $v\in V(G)$ is called a \it{central vertex} if $N^2(v)=\{v\}$.
 \item For every star $S_u\cup \{u\},$ its leaves (that is, vertices of $S_u$) will be called
the \textit{exterior  vertices} of this star.
\item If $P_k$ is a set of vertices, picked up at the $k^{th}$ stage, we will write $N_k$ to denote the set of pseudo-neighbors and pseudo-neighbors of pseudo-neighbors of the vertices of $P_k.$  That is,
    $$N_k:=\bigcup_{v\in P_k} \big(N(v)\cup N^2(v)\big).$$
\end{enumerate}\end{definition}
\begin{example}
If $P_k=\{a,b,c\}$ is a triplet, picked up at the $k^{th}$ stage, then
    $$N_k:=N(a)\cup N^2(a)\cup N(b)\cup N^2(b)\cup N(c)\cup N^2(c).$$
\end{example}
\begin{examples}
\begin{enumerate}
\item If $S_u$ is defined, then $u$ is a central vertex.
\item If $uv$ is an edge in the matching $M,$ then both $u\nd v$ are central vertices.
\item If $\{x,y,z\}$  is a triangle, as in Definition \ref{nndefn}, then $x$ is a central vertex.
\end{enumerate}
   \end{examples}

It must be noted that this notion of pseudo-neighbors is personal (dependent upon Algorithm \ref{alg41}) to this paper and should not be confused with the general neighborhood of a vertex of a graph.

 \begin{definition} Let $r>0$ be arbitrary and
$\del:= \frac{r}{6n}.$  For $v \in V(G),$ define
\begin{align*}
  m(v)&:=  min\big\{m:P_m\cap (N(v)\cup N^2(v))\ne \emptyset \big\}\\
&=min\big\{m:v \in P_m\mbox{ or }v_1\in P_m; v_1\in N(v)\mbox{ or }v_2\in P_m; v_2\in N^2(v)\big\}.
\end{align*}
Define $r(v):=r-2\delta m(v).$
 \end{definition}
 \begin{notations}
\begin{enumerate}
\item If $v\in V(G),$  then $nv$ will denote an arbitrary vertex in $N(v).$
\item If $v\in V(G),$ then $n^2v$ will denote an arbitrary vertex in $N^2(v).$
\item Let $P_k=\{p_k,q_k,s_k\}$ be a triplet, picked up at the $k^{th}$ stage, through the recurssion of step \ref{CharVerS24} of Algorithm \ref{alg41} such that $p_kq_k\in E(G).$
    Assume that there exist some $v\in V(G)$ \st $\{p_k,q_k,v\}\in \tcal.$
Then we will denote $v$ by $t_k.$
\end{enumerate}
 \end{notations}
\begin{note}\begin{myenu}
  \item For any $x\in S_u,$ we have $r(u)=r(x).$
\item For $\{x,y,z\}\in \T,$ we have $r(x)=r(y)=r(z).$
\item For $uv\in M,$ we have $r(u)=r(v).$
\end{myenu}\end{note}


\subsection{The Classification of Vertices of $G$}

Our  next step towards the progress on Conjecture \ref{conj} is to embed $G$ into a suitable Euclidean space.
We will establish   it in several different modules, for the convenience of readers, depending upon the classification of vertices picked up in Algorithm \ref{alg41}.

Let us revisit the
  classes of vertices of $G,$ picked up through  various steps of Algorithm \ref{alg41}.
  The  following table provides a brief summary for this purpose, mentioning
  the corresponding steps of  Algorithm \ref{alg41} through which each of these classes of vertices have been picked up.
This table also assigns a roman letter to each of these classes,  which will be used in next eight sections of this paper.
These roman letters will denote the corresponding modules to be explored in the sequel.

\begin{center}
\begin{tabular}{ |l|l|l| }
 \hline
Sr. No. & Class of vertices & Corresponding steps of Algorithm \ref{alg41}  \\
  \hline
I & Random pairs or singletons &  \ref{CharVerS18}-\ref{CharVerS19}, \ref{CharVerS30}, \ref{CharVerS37}   and  \ref{CharVerS45}  \\
II & The   residual sets & \ref{CharVerS27},  \ref{CharVerS32} or \ref{CharVerS38}   \\
III & Independent pairs $\{p,q\},$ where $ pq\notin E(G)$ & \ref{CharVerS40}   \\
IV & The triangles $\{p',q',s'\}$  &\ref{CharVerS43}   \\
V &  The case of $\{v_0,w_1,w_2\}$  & \ref{CharVerS33}   \\
VI & The case of  $\{w_0, v_1, v_2\}$  &  \ref{CharVerS35}   \\
VII & Triplets having exactly one edge & \ref{CharVerS24}  \\
VIII & The super triplets &   \ref{CharVerS07}-\ref{CharVerS17} and  \ref{CharVerS22}. \\
 \hline
\end{tabular}
\end{center}


By \textit{super triplets}, we refer to the   triplets $P_k,$
picked up  till step (\ref{CharVerS22}) of Algorithm \ref{alg41}.
\Nt each of these $P_k$ consists of independent vertices and  no two vertices out of these  belong to a single star.
Independent triplets are also picked up at step \ref{CharVerS33}, but there these triplets have    two vertices from a single star. That is why we need to coin separate terminology and notations to distinguish between these cases.

 \subsection{The Strategy}

Our strategy  to establish Theorem \ref{MainThmSIG3} can be briefly outlined as follows.

Corresponding to each of the picked sets of vertices $P_k,$
 we will suitably assign some Euclidean dimensions $c_k$ to all the vertices of $G,$
in order  to satisfy our requirements.
This will lead to an embedding of our graph into a suitable Euclidean space, as a sphere of influence graph.

Just for the sake of convenience, we will use the same symbol $v$ for the image of $v,$ under this embedding.
In this sense, the vertex set $V(G)$ will be  projected in an Euclidean space endowed with sup metric.
We will  prove that the sphere of influence graph of this mapped vertex set is isomorphic to our given graph $G.$

We will use the notation $|c_k(u)-c_k(v)|; u, v\in G,$ even when $c_k$ represents multiple  Euclidean dimensions. In that case, as an abuse of notation,  it will represent the sup-norm in those  dimensions given by $c_k.$
Whenever we require to deal with individual dimensions of $c_k(v),$ we will denote them separately by $c_k[0](v), c_k[1](v) \dots.$

Note that given any category of the set $P_k$
of vertices picked up through Algorithm \ref{alg41},
there may exist a graph $G$
such that Algorithm \ref{alg41} provides only that category of picked up vertices only.


We provide a set of inequalities, to be established in next eight sections.
In each of the following, let $P_k$ be any set of vertices, picked up  through Algorithm \ref{alg41}.

\noindent If $u\in V(G) \nd nu\in N(u),$ then
\begin{equation}\label{LemEq1}
 |c_k(u)-c_k(nu)|  \leq r(u).
\end{equation}
If  $u\in P_k$  and  $v\in \displaystyle\bigcup_{l\leq k} P_l\sm \{u\},$ then
\begin{equation}\label{LemEq2}
 |c_k(v)-c_k(u)|  \geq \max\{r(u),r(v)\}.
\end{equation}
\noindent If  $u\in P_k\nd  v\in \displaystyle\bigcup_{l\leq k} P_l \sm \{u\}$  such that $uv\not\in E(G)$   and $\{u,v\}\subseteq S_a \fos a\in V(G),$  then
\begin{equation}\label{LemEq3}
|c_k(v)-c_k(u)|  \geq r(v)+r(u).
\end{equation}
If  $u\in P_k$  and  $v\in \displaystyle\bigcup_{l\geq k} P_l\sm \{u\} $ such that $uv\not\in E(G)$   and $\{u,v\}\not \subseteq S_w \foa w\in V(G),$  then
\begin{equation}\label{LemEq4}
|c_k(v)-c_k(u)|  \geq r(v)+r(u).
\end{equation}
If   $u\in V(G)\nd v\in V(G) \sm \{u\}$ such that $uv\in E(G), $  then
\begin{equation}\label{LemEq5}
|c_k(v)-c_k(u)|  < r(v)+r(u).
\end{equation}

\noindent First of all we note that   the assertions (\ref{LemEq3})   and  (\ref{LemEq4}), for each $k,$
will establish the following.
\begin{equation*}\label{LemEq7}
|c_k(v)-c_k(u)|  \geq r(v)+r(u)\foa u, v\in V(G)\mbox{ such that }uv\not\in E(G).
\end{equation*}

\noindent   Further, note that inequality  (\ref{LemEq2}) for all $k,$ implies the following
\begin{equation*}
 \rho(u,v) \geq \max\{r(u),r(v)\}\foa u,v\in V(G), u\ne v.
\end{equation*}


\noindent Hence, proving (\ref{LemEq1})   and  (\ref{LemEq2}) for each $k$ will lead to
$$r_u=r(u)\foa u\in V(G).$$
\Csq  the   assertions  (\ref{LemEq1})   to  (\ref{LemEq5})   will jointly conclude the following
$$v_1v_2\in E(G)\mbox{  if and only if }\rho(v_1, v_2) < r_{v_1}+r_{v_2}\foa v_1, v_2\in V(G).$$
This will establish that our embedding of $G$ is a sphere of influence graph on a set of vertices in some Euclidean space.
Finally, in last section of this paper, we will count the dimensions involved in this procedure.
An eager reader can have a look at Theorem  \ref{MainThmSIG3}.

In each of the following modules, we will be establishing these assertions  (\ref{LemEq1})   to  (\ref{LemEq5}), for specified sets of vertices picked up through Algorithm \ref{alg41}. Every module will start with a reference to a class of vertices, as per the previous table. Then we will assign some suitable dimensions $c_k$ for each such $P_k.$ After that
we divide each module into  five subsections, to prove each of the assertions
(\ref{LemEq1})   to  (\ref{LemEq5}).
This will provide a justification of our dimension assignment, for each set of vertices $P_k.$
As remarked earlier, it will conclude the required embedding.

For the sake of convenience of the readers, we will be recalling these inequalities before starting their proofs.

It must be noted that each of  these modules can be verified in any order.
We will be discussing  these modules, in the  increasing level of the corresponding technicalities involved.

For example, the super triplets and the  triplets having exactly one edge
appear earlier than most of the other cases, in the  Algorithm \ref{alg41}.
However, we will be dealing with these cases in the last modules, as these involve some of the most technical manipulations.

Due to abundance of notations and technical jargon, we will be writing each of (\ref{LemEq1}) to (\ref{LemEq5}), in terms of lemmas in the last two modules.

We wind up this section, with a few little notes which would be used again and again in the subsequent modules.

If $u\in P_i,$ then our definitions and notations imply that $$m(u)\leq i\mbox{ and thence }r(u)\geq r-2\del i.$$
For each $u\in V(G),$ the set $N(u)\cup N^2(u)$ is closed under the pseudo-neighborhood operation.
Therefore, for (\ref{LemEq1}), it is enough to establish
$$|c_k(nu)-c_k(n^2u)|\leq r(u)\foa  nu\in N(u)\ndfoa n^2u\in N^2(u).$$
Hence, if (\ref{LemEq1}) is satisfied for all  $u\in N(w),$ then (\ref{LemEq1}) also holds for
for all $u\in N^2(w).$

%

 \subsection{Sample Calculations}

In this subsection, we outline a few sample cases for the inequalities to be verified, at different stages, in the next eight sections of this paper. Our aim is to give the reader a gist of our upcoming  calculations.
It is a very little sample and  leads to most of the later inequalities.


\begin{myenu}
\item $\delta  \leq \frac{r}{6}.$
\begin{proof}
   Since our graph $G$ has no isolated vertices,   $n\geq 1.$
Hence $\delta =\frac{r}{6n}\leq \frac{r}{6}.$
\end{proof}
\item  If $P_k$ exists and    $v\in N_k,$ then $r(v)>r-2\delta (k+1).$
\begin{proof}
Since  $v\in N_k,$ \wv $m(v)\leq k <k+1.$
Hence  $r(v)=r-2\delta m(v)>r-2\delta (k+1).$
\end{proof}
\item  For all $v\in V(G),$ \wv $\frac{2n}{3}<r(v)\leq r.$
\begin{proof}
Recall that    $r(v)=r-2\delta m(v)$ and $\delta=r/6n.$
Since  $ m(v)\geq 0,$ \wv  $r(v)\leq r.$ Since  $ m(v)<n,$ \wv
\begin{equation*}
   r(v)=r-\frac{r}{3n}m(v)> r-\frac{r}{3}=\frac{2r}{3}.\qedhere
   \end{equation*}
\end{proof}
\item    If $P_k$ exists, then $\frac{2r}{3}<r-2\delta (k+1)<r.$
\begin{proof}
 Since $P_k$ exists, our construction implies that $k+1<n.$
Hence
$$r-2\delta (k+1)= r- \frac{r}{3n}(k+1)>r- \frac{r}{3n}n=\frac{2r}{3}.$$
The other inequality is immediate, as $k+1>0.$
%
%
 \end{proof}

\item  Let $u,v\in V(G)$ \st $u\in P_k.$ Let $c_k(v)=(-r(v),-r(v))$ and
$$c_k(u)=(-r+2\delta (k+1)-r(u),-r+2\delta (k+1)+r(u)-\delta ).$$
Then
\begin{equation}\label{harnoortangkarda}
 |c_k(u)-c_k(v)|\geq \max\{r(u),r(v)\}.
 \end{equation}
 (This will be used in Case 2.1 on page  \pageref{p44ineq5}.)
\begin{proof}
\label{pageconj14}
Since $u\in P_k,$ \wv    $m(u)\leq k.$ Hence $r(u)\geq r-2\delta k.$\Nt
$$  |c_k(u)-c_k(v)|= \max\big\{r(u)+r(v)-r+2\delta (k+1)-\delta, r(u)+r-2\delta (k+1)-r(v) \big\}.$$
Consider the following cases.
\begin{myenu}[label=(\roman*)]
\item \bm{$r(v)\leq r-2\delta (k+1).$}
Then
$r(u)+r-2\delta (k+1)-r(v)\geq r(u) .$
Already \wv $r(u)\geq r-2\delta k,$ which implies  $r(u)\geq r(v).$ Hence
$$r(u)+r-2\delta (k+1)-r(v)\geq r(u)\geq r(v).$$
\aThr $r(u)+r-2\delta (k+1)-r(v)\geq \max\{r(u),r(v)\}.$
\item \bm{$r(v)>r-2\delta (k+1).$} Then  $r(v)\geq r-2\delta k.$
Hence
$$r(u)+r(v)-r+2\delta (k+1)-\delta
\geq r(u)+\delta > r(u).$$
Similarly, as we already have  $r(u)\geq r-2\delta k,$ analogously \wo
$$r(u)+r(v)-r+2\delta (k+1)-\delta> r(v).$$
Hence, \wo
$$r(u)+r(v)-r+2\delta (k+1)-\delta >\max\{r(u),r(v)\}.$$
\end{myenu}
This proves  (\ref{harnoortangkarda}), in both of the cases. \qedhere

\end{proof}
\end{myenu}

These are just few sample calculations. Anyone interested in more such details is referred to  the last module in section \ref{sectionM8},  where we have provided detailed proofs of several other such assertions.

However the authors are of the opinion that including all such calculations would unnecessarily divert the main theme of the paper. Therefore we are avoiding these kind of routine manipulations in the first seven modules, most of which simply follow by the definitions of $r(u)$ and $m(u)$ etc.

\section{\bf Module I: Random Pairs or Singletons}\label{subcaseIIa}

Consider the sets of vertices picked up through the recursions of   any of the steps
 \ref{CharVerS18}-\ref{CharVerS19}, \ref{CharVerS30},  \ref{CharVerS37}  or \ref{CharVerS45} of   Algorithm \ref{alg41}.
In this case, we reserve one dimension for each vertex in such a set.
Let $P_k$ be a pair or a singleton set, picked through any of these steps.

If $P_k$ is the singleton set, $\{w\},$ we define $c_k$ on the vertices of $G,$ as follows.
\begin{myenu}[label=Step \arabic*.]
\item   Define   $c_k(w):=-r(w).$
\item   Define  {$c_k$ of  any $nw\in N(w)$} as follows: $c_k(nw):=0.$
\item    Define {$c_k$ of   $v\in V(G)\setminus \{w,N(w)\}$} as follows.
$$c_k(v):=\left\{\begin{array}{ll}
r(v)-\delta 		& \mbox{; $vw\in E(G),$ }\\
r(v)				& \mbox{; $vw\notin E(G).$ }
\end{array}\right.$$
\end{myenu}
If $P_k$ is a pair $\{p,q\},$ then for each of $p\nd q,$ we reserve individual separate dimensions, as above.
We denote these two dimensions by $c_k[0]$  and $c_k[1].$

Next we justify this choice of new dimensions. Note that it is enough to establish that only for the case of singleton sets $P_k.$ Hence, without loss of generality, we assume that $P_k=\{w\}.$

\subsection{Proof of (\ref{LemEq1})}
Let $v\in V(G) \nd nv\in N(v).$ We need to establish that
\begin{equation*}
 |c_k(v)-c_k(nv)|  \leq r(v).
\end{equation*}
This is obvious, due to our choice of $\del\nd r(v),$ as we have only the following cases.
\begin{myenu}[label=Case \arabic*.]
\item \bm{$v\in N_k=N(w)\cup N^2(w).$}
Then one of $c_k(v)\nd c_k(nv)$ is zero and the other belong to the set  $\{-r(v),  r(v), r(v)-\del\}.$
 \item \bm{$v\in V(G)\sm N_k.$} Then both $c_k(v)\nd c_k(nv)$ belong to the set  $\{r(v), r(v)-\del\}.$
\end{myenu}

\subsection{Proof of (\ref{LemEq2})}

Let  $u\in P_k$  and  $v\in \displaystyle\bigcup_{l\leq k} P_l\sm \{u\}.$
 We need to  prove that
\begin{equation*}
 |c_k(v)-c_k(u)|  \geq \max\{r(u),r(v)\}.
\end{equation*}
Since $u\in P_k=\{w\},$ we have $u=w$ and therefore $c_k(u)=-r(u).$
Further, we note that $c_k(v)$ belongs to $\{0, r(v), r(v)-\del\}.$


In case $c_k(v)\in \{  r(v), r(v)-\del\},$ the above inequality is satisfied due to our definitions of $r(u), r(v)\nd \del.$ If $c_k(v)=0,$ then $v\in N(w)=N(u)\wimp r(u)=r(v).$ This establishes the above inequality and hence (\ref{LemEq2}).

\subsection{Proof of (\ref{LemEq3})}

Let  $u\in P_k\nd  v\in \displaystyle\bigcup_{l\leq k} P_l \sm \{u\}$  such that $uv\not\in E(G)$   and $\{u,v\}\subseteq S_a \fos a\in V(G).$
 We need to  prove that
$$|c_k(v)-c_k(u)|  \geq r(v)+r(u).$$

%
Since $u\in P_k=\{w\},$ \wv $c_k(u)=-r(u).$
Further $uv\notin E(G)$ implies that $c_k(v)=r(v).$


\subsection{Proof of (\ref{LemEq4})}

Let  $u\in P_k$  and  $v\in \displaystyle\bigcup_{l\geq k} P_l\sm \{u\} $ be such that $uv\not\in E(G)$   and $\{u,v\}\not \subseteq S_w \foa w\in V(G).$
We need to  prove that
$$|c_k(v)-c_k(u)|  \geq r(v)+r(u).$$

Again, as above, we have $c_k(u)=-r(u)\nd c_k(v)=r(v).$

\subsection{Proof of (\ref{LemEq5})}
Let   $u\in V(G)\nd v\in V(G) \sm \{u\}$ such that $uv\in E(G).$  It is enough to prove that
$$|c_k(v)-c_k(u)|  < r(v)+r(u).$$
 This inequality is clearly satisfied as  the following are the only possibilities.
 \begin{myenu}[label=Case \arabic*.]
\item \bm{$u=w.$}
Then $c_k(u)=-r(u)\nd c_k(v)\in \{0,r(v)-\delta \}.$
\item \bm{$u\neq w.$} Then
 $c_k(u)\in \{0,r(u)-\delta ,r(u)\}\nd c_k(v)\in \{0,r(v)-\delta ,r(v)\}.$
\end{myenu}

\section{\bf Module II: The Residual Sets}\label{subcaseIIa}
Consider the sets of vertices picked up through the recursion of   any of the steps  \ref{CharVerS27}, \ref{CharVerS32} or \ref{CharVerS38} of Algorithm \ref{alg41}. These will be termed as the residual sets, throughout this module.
Note that Algorithm \ref{alg41} will yield at most one of these residual sets.
Let $P_k$ denote that set of vertices, say $P_k=\{e_1,\dots, e_m\}.$
We reserve $\lceil log_2 (m+1)\rceil$ dimensions for $P_k.$

First we present a  result, due to our construction, which will be used frequently in this paper.

\begin{prop}\label{propII}
Let $P_k$ be a residual set, $u\in P_k$ and $v\in \bigcup_{l>k} P_{l}.$ Then $uv\in E(G).$
\end{prop}
\begin{proof}
There are three cases to consider.
\begin{myenu}[label=Case \arabic*.]
\item \bm{$P_k$}\textbf{ is picked up through step \ref{CharVerS27} of Algorithm \ref{alg41}.}

\Nt in this case,  \ref{alg41} ensures that $|P_k|\geq 2.$
Let $u_1, u_2$ be any two vertices from $P_k,$ which belong to different stars and  $v\in  \cup_{l>k} P_{l}.$
Note that all the possible vertices of triplets  with no edge or exactly one edge had already been picked up prior to $P_k,$
using steps  \ref{CharVerS22} and \ref{CharVerS24} of  Algorithm \ref{alg41}.
Further $u_1u_2\notin E(G),$ being  exterior vertices of our stars.
Hence,  we conclude that the vertices $u_1,u_2,v$ have exactly two edges, $u_1v\nd u_2v.$
Since $u_1$ and $u_2$ were arbitrary vertices from $P_k,$ belonging to different stars, we conclude the result for this case.
\item \bm{$P_k$}\textbf{ is picked up through step \ref{CharVerS32} of Algorithm \ref{alg41}.}

In this case, no set of vertices is picked after $P_k$ and therefore the result is vacuously true.
\item \bm{$P_k$}\textbf{ is picked up through step \ref{CharVerS38} of Algorithm \ref{alg41}.}

If there is no $P_l$ with $l>k,$ the  result is vacuously true.
Otherwise, due to  step \ref{CharVerS37} of Algorithm \ref{alg41}, there are at least two elements in $\cup_{l>k} P_{l}.$
Let $u\in P_k$ and  $v_1,v_2\in  \cup_{l>k} P_{l}.$
Due to steps \ref{CharVerS22} and  \ref{CharVerS24}, $\{u, v_1,v_2\}$ either has exactly two edges or exactly three edges.
If it has exactly three edges, then we have nothing to prove.
If $uv_1\notin E(G),$ then again we have a contradiction due to our choice in step \ref{CharVerS35}.
Hence, $uv_1\in E(G).$ Similarly, \wo    $uv_2\in E(G).$
%

\end{myenu}
This  proves the result in all cases.
\end{proof}

 Now, we assign suitable dimensions to the vertices of $G,$ corresponding to  $P_k=\{e_1,\dots, e_m\}.$
Let $v\in V(G).$  We define $c_k(v)$ as a $\lceil log_2 (m+1)\rceil$-dimensional point.

Note that there exist at least $m+1$ distinct $\lceil log_2 (m+1)\rceil$-dimensional points, with each co-ordinate either $1\mbox{ or }-1.$ Let $f_1,f_2,\dots,f_{m+1}$ be any such distinct points.
We define $c'(v)$ as a $\lceil log_2 (m+1)\rceil$-dimensional point as follows.
$$c'(v):=\left\{\begin{array}{ll}
f_{i} & \mbox{; $v=e_{i},$ for some $1\leq i\leq m,$ }\\
f_{m+1} & \mbox{; $v\in V(G)\setminus P_k.$}\end{array}\right.$$
Finally, we define
$$c_k(v):=\left\{\begin{array}{ll}
r(v)c'(v) & \mbox{; $v\in P_k,$ }\\
(0,0,\dots,0) & \mbox{; $v\in N(P_k),$ }\\
r(v)c'(v) & \mbox{; $v\in N^2(P_k)\setminus P_k,$}\\
(r(v)-\delta )c'(v) & \mbox{; otherwise.}\end{array}\right.$$

\subsection{Proof of (\ref{LemEq1})}
Let $u\in V(G) \nd nu\in N(u).$ We need to establish that
\begin{equation*}
 |c_k(u)-c_k(nu)|  \leq r(u).
\end{equation*}
\Nt we have the following cases.
\begin{myenu}[label=Case \arabic*.]
\item \bm{$u\in N_k.$} Then $nu\in V_k$ and hence both $c_k(u)\nd c_k(nu)$ belong to the set
$$\big\{(0,0,\dots,0),r(u)c'(u)\big\}.$$
\item \bm{$u\in V(G)\sm N_k.$} Then
\begin{align*}
   c_k(v)=&(r(v)-\delta )c'(v)=(r(v)-\delta )f_{m+1}\\
   c_k(nv)=&(r(nv)-\delta )c'(nv)=(r(v)-\delta )f_{m+1}.
\end{align*}
\end{myenu}
Due to our choice of $c'(u),$  the above inequality is trivially satisfied.

\subsection{Proof of (\ref{LemEq2})}
Let  $u\in P_k$  and  $v\in \displaystyle\bigcup_{l\leq k} P_l\sm \{u\}.$
 We need to  prove that
\begin{equation*}
 |c_k(v)-c_k(u)|  \geq \max\{r(u),r(v)\}.
\end{equation*}
Since $u\in P_k,$ \wv $c_k(u)=r(u)c'(u).$
We have the following cases.
\begin{myenu}[label=Case \arabic*.]
\item \bm{$v\in N_k.$} Then  $c_k(v)\in\{r(v)c'(v),(0,0,\dots,0) \}=\{r(u)c'(v),(0,0,\dots,0)\}.$
\item \bm{$v\in V(G)\sm N_k.$} Then  $c_k(v)=(r(v)-\del)c'(v).$
\end{myenu}
Our choice of $c'(u)$ ensures that the above inequality is   satisfied.

%
%

\subsection{Proof of (\ref{LemEq3})}

Let  $u\in P_k\nd  v\in \displaystyle\bigcup_{l\leq k} P_l \sm \{u\}$  such that $uv\not\in E(G)$   and $\{u,v\}\subseteq S_a \fos a\in V(G).$
 We need to  prove that
$$|c_k(v)-c_k(u)|  \geq r(v)+r(u).$$
Then $u,v \in N(a)$ and hence $v\in N^2(u).$ \Csq
$c_k(u)=r(u)c'(u)\nd c_k(v)= r(v)c'(v).$
Therefore the above inequality is satisfied.

\subsection{Proof of (\ref{LemEq4})}

Let  $u\in P_k$  and  $v\in \displaystyle\bigcup_{l\geq k} P_l\sm \{u\} $ be such that $uv\not\in E(G)$   and $\{u,v\}\not \subseteq S_w \foa w\in V(G).$\\
We need to  prove that
$$|c_k(v)-c_k(u)|  \geq r(v)+r(u).$$
As above, by Proposition \ref{propII}, we obtain $v\in P_k.$
Hence
$c_k(u)=r(u)c'(u)\nd c_k(v)= r(v)c'(v).$
Therefore the above inequality is satisfied.

%
%
%

\subsection{Proof of (\ref{LemEq5})}

Let   $u\in V(G)\nd v\in V(G) \sm \{u\}$ such that $uv\in E(G).$  It is enough to prove that
$$|c_k(v)-c_k(u)|  < r(v)+r(u).$$
Assume that $u\in P_i\nd v\in P_j.$
We have the following cases.
\begin{myenu}[label=Case \arabic*.]
\item \bm{$k\notin \{i,j\}.$} Then $c'(u)=c'(v)=f_{m+1}.$ Let
$$A(u):=\{(0,0,\dots,0),r(u)c'(u),(r(u)-\delta )c'(u)\}
.$$
Then  $c_k(u)\in A(u)\nd c_k(v)\in A(v).$
\item \bm{ $k\in \{i,j\}.$} Without loss of generality, we assume that $k=i.$
Then $$c_k(u)=r(u)c'(u)\nd c_k(v)\in \{(0,0,\dots,0),(r(v)-\delta )c'(v)\}.$$
\end{myenu}
Due to our choice of $c'(u),$  the above inequality is  satisfied.

\section{\bf Module III: $\{p,q\},$ where $ pq\notin E(G)$}

Consider the sets of pairs $\{p,q\}$ \st $pq\notin E(G),$  picked up through the recursion of    step \ref{CharVerS40} of   Algorithm \ref{alg41}. In this case, we reserve one dimension for each such pair.
Let $P_k=\{p_k,q_k\}$ be any such pair.

Analogous to Proposition \ref{propII}, for this $P_k,$ we have  the following result.
\begin{prop}\label{propIII}
   If $u\in P_k\nd v\in \bigcup_{l>k} P_{l},$ then $uv\in E(G).$
   \end{prop}
\begin{proof}
Write $P_k=\{p_k,q_k\}.$
Note that all the possible vertices of triplets  with no edge or exactly one edge had already been picked up prior to $P_k,$
 using steps  \ref{CharVerS22} and \ref{CharVerS24} of  Algorithm \ref{alg41}.
Further $p_kq_k\notin E(G).$
Hence,  we conclude that the vertices $p_k,q_k,v$ have exactly two edges, $p_kv\nd q_kv.$
This proves the result.
\end{proof}
Below we assign a corresponding   dimension $c_k$ on the vertices of $G,$ as follows.

\begin{myenu}[label=Step \arabic*.]
\item Define $c_k(p_k):=-r(p_k)-\frac{r}{2}+\delta (k+1).$
\item Define  $c_k(q_k):=r(q_k)+\frac{r}{2}-\delta (k+1).$
\item If $np_k\in N(p_k),$ define
$c_k(np_k):=-\frac{r}{2}+\delta (k+1).$
\item If $n^2p_k\in N^2(p_k)\setminus \{p_k\},$ define
$c_k(n^2p_k):=-\frac{r}{2}+\delta (k+1).$
\item If $nq_k\in N(q_k),$ define $c_k(nq_k):=\frac{r}{2}-\delta (k+1).$
\item If $n^2q_k\in N^2(q_k)\setminus \{q_k\},$ define
$c_k(n^2q_k):=\frac{r}{2}-\delta (k+1).$
\item Finally, if $v\in V(G)\setminus V_k,$ define
$c_k(v):=0.$
\end{myenu}

\subsection{Proof of (\ref{LemEq1})}
Let $u\in V(G) \nd nu\in N(u).$ We need to establish that
\begin{equation*}
 |c_k(u)-c_k(nu)|  \leq r(u).
\end{equation*}
\Nt we have the following cases.
\begin{myenu}[label=Case \arabic*.]
\item \bm{$u \in N(p_k).$} Then $nu \in N^2(p_k)$ and we have
$$c_k(u)=-\frac{r}{2}+\delta (k+1)\nd c_k(nu)\in \bigg\{-\frac{r}{2}+\delta (k+1),-r(p_k)-\frac{r}{2}+\delta (k+1)\bigg\}.$$
Therefore, we have
$|c_k(u)-c_k(nu)|$ is either $0\mbox{ or }|-r(p_k)|,$ which is equal to $r(u).$ Hence the required  inequality is satisfied.
\item \bm{$u \in N^2(p_k).$} Then $nu \in N(p_k).$ This case is equivalent to the above case.
\item \bm{$u \in N(q_k)\cup N^2(q_k).$} Similar to the above cases.
\item \bm{$u\in V(G)\sm  V_k.$}\ Then $nu \in V(G)\sm  V_k$ and hence $c_k(u)=0=c_k(nu).$
\end{myenu}

\subsection{Proof of (\ref{LemEq2})}

Let  $u\in P_k$  and  $v\in \displaystyle\bigcup_{l\leq k} P_l\sm \{u\}.$
 We need to  prove that
\begin{equation*}
 |c_k(v)-c_k(u)|  \geq \max\{r(u),r(v)\}.
\end{equation*}

\begin{myenu}[label=Case \arabic*.]
\item \bm{$u=p_k.$} Then $c_k(u)=-r(u)-\frac{r}{2}+\delta (k+1)\nd $
$$c_k(v)\in \bigg\{r(v)+\frac{r}{2}-\delta (k+1),-\frac{r}{2}+\delta (k+1), \frac{r}{2}-\delta (k+1),0\bigg\}.$$
\item \bm{$u=q_k.$} Similar to the above case.
\end{myenu}

\subsection{Proof of (\ref{LemEq3})}

Let $u\in P_k\nd  v\in \displaystyle\bigcup_{l\leq k} P_l \sm \{u\}$ be such that $uv\not\in E(G)$   and $\{u,v\}\subseteq S_w \fos w\in V(G).$
We need to  prove that
$$|c_k(v)-c_k(u)|  \geq r(v)+r(u).$$

Since $u\in P_k,$ which has been chosen at step
\ref{CharVerS40} of   Algorithm \ref{alg41}, it is not an exterior vertex of any star.
Hence   there are no such vertices $u\nd v$ and consequently the above inequality is vacuously true.

\subsection{Proof of (\ref{LemEq4})}

Let $u\in P_k$  and  $v\in \displaystyle\bigcup_{l\geq k} P_l\sm \{u\} $ be such that $uv\not\in E(G)$   and $\{u,v\}\not \subseteq S_w \foa w\in V(G),$  We need to  prove that
$$|c_k(v)-c_k(u)|  \geq r(v)+r(u).$$

Since $uv\notin E(G),$ applying Proposition \ref{propIII}, \wo  $u, v\in P_k.$ Without loss of generality, assume that $u=p_k\nd v=q_k.$ Then
$$c_k(u)= -r(u)-\frac{r}{2}+\delta (k+1)\nd c_k(v)=r(v)+\frac{r}{2}-\delta (k+1).$$
Hence  inequality (\ref{LemEq4}) is valid for this module.

\subsection{Proof of (\ref{LemEq5})}

Let   $u\in V(G)\nd v\in V(G) \sm \{u\}$ such that $uv\in E(G).$  We need to prove that
$$|c_k(v)-c_k(u)|  < r(v)+r(u).$$

%
%
Consider the set
$$A(u):=\bigg\{-r(u)-\frac{r}{2}+\delta (k+1),r(u)+\frac{r}{2}-\delta (k+1),-\frac{r}{2}+\delta (k+1), \frac{r}{2}-\delta (k+1),0\bigg\}.$$
Then $c_k(u)\in A(u)\nd c_k(v)\in A(v).$ Also, observe that
$$|c_k(u)-c_k(v)| > r(u)+r(v)\mbox{ only when }\{u,v\}=\{p_k,q_k\}.$$
In these cases, we have  $uv\notin E(G),$ a contradiction.
Hence
(\ref{LemEq5}) is satisfied.

%
%
%
%


\section{\bf Module IV: The triangles $\{p',q',s'\}$}

Consider the sets of triangles picked up through the recursion of    step \ref{CharVerS43} of   Algorithm \ref{alg41}.
We reserve two dimensions for each such triangle.
Let $P_k=\{p'_k,q'_k,s'_k\}$ be any such triangle.

Below, we assign the corresponding  dimensions $c_k$ on the vertices of $G.$

\begin{myenu}[label=Case \arabic*.]
\item \textit{No vertex of $P_k$ is a pseudo-neighbor or pseudo-neighbor of pseudo-neighbor  of any other vertex of $P_k.$}

In this case, we write $P_k$ as $\{x,y,z\},$ where $$r(x):=\max\{r(p'_k),r(q'_k),r(s'_k)\}.$$
If multiple vertices $v$ of $P_k$ have maximum $r(v),$ then call any of these    as $x.$
Let $y\nd z$ denote the other two vertices of $P_k.$ \Nt we have
$$x\notin N(y)\cup N^2(y)\cup N(z)\cup N^2(z)\mbox{ and }y\notin  N(z)\cup N^2(z).$$
For $v\in V(G),$ we define
 $$c_k(v):=\left\{\begin{array}{ll}
(0,0) & \mbox{; $v=x,$ }\\
(0,r) & \mbox{; $v=y,$ }\\
(r,0) & \mbox{; $v=z,$ }\\
(r(x),r(x)) & \mbox{; $v\in N(x),$ }\\
(r(y),r) & \mbox{; $v\in N(y),$ }\\
(r,r(z)) & \mbox{; $v\in N(z),$ }\\
(r,r) & \mbox{; otherwise}. \end{array}\right.$$
\item \textit{One vertex of $P_k$ is a pseudo-neighbor or pseudo-neighbor  of pseudo-neighbor  of another vertex of $P_k,$ but the}
\textit{third vertex is neither a pseudo-neighbor nor  a pseudo-neighbor of pseudo-neighbor  of the other two vertices.}

Without loss of generality, assume that
$$p'_k\in N(q'_k)\cup N^2(q'_k)\mbox{ and }s'_k\notin N(p'_k)\cup N^2(p'_k)\cup N(q'_k)\cup N^2(q'_k).$$
 For $v\in V(G),$  we define
$$c_k(v):=\left\{\begin{array}{ll}
(0,r-r(v)) & \mbox{; $v=p'_k,$ }\\
(0,r) & \mbox{; $v=q'_k,$ }\\
(r(v),r) & \mbox{; $v\in N(p'_k)\setminus \{q'_k\},$ }\\
(r(v),r) & \mbox{; $v\in N(q'_k)\setminus \{p'_k\},$ }\\
(r,0) & \mbox{; $v=s'_k,$}\\
(r,r(v)) & \mbox{; $v\in N(s'_k),$}\\
(r,r) & \mbox{; otherwise}. \end{array}\right.$$

 \item \textit{$\{p'_k,q'_k,s'_k\}\in \tcal.$} Define
$$c_k(v):=\left\{\begin{array}{ll}
(-r(v),0) & \mbox{; $v=p'_k,$ }\\
(-r(v),-r(v)) & \mbox{; $v=q'_k,$ }\\
(0,-r(v)) & \mbox{; $v=s'_k,$ }\\
(r(v)-\delta ,r(v)-\delta ) & \mbox{; otherwise}. \end{array}\right.$$

\end{myenu}

Next, we will prove each of  the assertions (\ref{LemEq1}) to  (\ref{LemEq5}), for this $P_k.$
In each  subsection of this module,   we will be dealing with the above cases separately.
However, we won't be repeating all the details of these cases therein.

\subsection{Proof of (\ref{LemEq1})}
Let $v\in V(G) \nd nv\in N(v).$ We need to establish that
\begin{equation*}
 |c_k(v)-c_k(nv)|  \leq r(v).
\end{equation*}

%
We now analyze this case by case as follows.

\begin{myenu}[label=Case \arabic*.]
\item \bm{$P_k=\{x,y,z\}.$}
\begin{myenu}
\item\label{IVcasdaeeq1a1} \bm{$u \in N(x).$} Then $nu \in N^2(x),
c_k(u)=(r(u),r(u))\nd c_k(nu)\in \{(0,0),(r,r)\}$
\item\label{IVcasdaeeq1a2} \bm{$u \in N(y).$} Then $nu \in N^2(y),$
$c_k(u)=(r(u),r)\nd c_k(nu)\in \{(0,r),(r,r)\}$
\item\label{IVcasdaeeq1a3} \bm{$u \in N(z).$} Then $nu \in N^2(z),$
$c_k(u)=(r,r(u))\nd c_k(nu)\in \{(r,0),(r,r)\}$
\item \bm{$u \in N^2(x).$} Then $nu \in N(x).$ This case follows from  \ref{IVcasdaeeq1a1}
\item \bm{$u \in N^2(y).$} This case follows from \ref{IVcasdaeeq1a2}
\item \bm{$u \in N^2(z).$} This case follows from \ref{IVcasdaeeq1a3}
\item \bm{$u\in V(G)\sm V_k.$}
Then $nu \in V(G)\sm V_k$
and we have $c_k(u)=(r,r)=c_k(nu).$
\end{myenu}
\item \bm{$p'_k\in N(q'_k)\cup N^2(q'_k)\mbox{ and }s'_k\notin N(p'_k)\cup N^2(p'_k)\cup N(q'_k)\cup N^2(q'_k).$}
\begin{myenu}
\item\label{IVcasdaeeq1a1s2e1} \bm{$u \in N(p'_k).$}
Then $nu \in N^2(p'_k)$ and we have
$$c_k(u)\in \{(0,r),(r(u),r)\}\nd c_k(nu)\in \{(0,r-r(u)),(0,r),(r(u),r)\}.$$
\item\label{IVcasdaeeq1a1s2e2} \bm{$u \in N(q'_k).$}
Then $nu \in N^2(q'_k)$ and we have
$$c_k(u)\in \{(0,r-r(u)),(r(u),r)\}\nd c_k(nu)\in \{(0,r-r(u)),(0,r),(r(u),r)\}.$$
\item\label{IVcasdaeeq1a1s2e3} \bm{$u \in N(s'_k).$}
Then $nu \in N^2(s'_k)$ and we have
$$c_k(u)=(r,r(u))\nd c_k(nu)\in \{(r,0),(r,r)\}.$$
\item \bm{$u \in N^2(p'_k).$}  Then $nu \in N(p'_k).$ This case follows from  \ref{IVcasdaeeq1a1s2e1}
\item \bm{$u \in N^2(q'_k).$}  Then $nu \in N(q'_k).$ This case follows from  \ref{IVcasdaeeq1a1s2e2}
\item \bm{$u \in N^2(s'_k).$}  Then $nu \in N(s'_k).$ This case follows from  \ref{IVcasdaeeq1a1s2e3}
\item \bm{$u\in V(G)\sm V_k.$} Then $nu\in V(G)\sm V_k$ and we have
$$c_k(u)=(r,r)\nd c_k(nu)=(r,r).$$
\end{myenu}
\item   \bm{$\{p'_k,q'_k,s'_k\}\in \tcal.$}
\begin{myenu}
\item \bm{$u\in \{p'_k,q'_k,s'_k\}.$} Then $nu\in \{p'_k,q'_k,s'_k\}.$
In this case
\begin{align*}
  c_k(u)&\in \{(-r(u),0),(-r(u),-r(u)),(0,-r(u))\}\\
\nd
c_k(nu)&\in \{(-r(u),0),(-r(u),-r(u)),(0,-r(u))\}.
\end{align*}
\item \bm{$u\in V(G)\sm \{p'_k,q'_k,s'_k\}.$} Then $ nu\in V(G)\sm \{p'_k,q'_k,s'_k\}$ and we have
$$c_k(u)=(r(u)-\delta ,r(u)-\delta )\nd c_k(nu)=(r(u)-\delta ,r(u)-\delta ).$$
\end{myenu}
\end{myenu}
%

\subsection{Proof of (\ref{LemEq2})}

Let  $u\in P_k$  and  $v\in \displaystyle\bigcup_{l\leq k} P_l\sm \{u\}.$
 We need to  prove that
\begin{equation*}
 |c_k(v)-c_k(u)|  \geq \max\{r(u),r(v)\}.
\end{equation*}
 We now study this case by case as follows.

  \begin{myenu}[label=Case \arabic*.]
\item \bm{$P_k=\{x,y,z\}.$}
 \begin{myenu}[label*= \arabic*.]
 \item If $u=x,$ then
 $c_k(u)=(0,0)$ and $$c_k(v)\in \{(0,r),(r,0),(r(v),r(v)),(r(v),r),(r,r(v)),(r,r)\}.$$
\item  If $u=y,$ then
 $c_k(u)=(0,r)$ and $$c_k(v)\in \{(0,0),(r,0),(r(v),r(v)),(r(v),r),(r,r(v)),(r,r)\}.$$
\item if $u=z,$ then
$c_k(u)=(r,0)$ and $$c_k(v)\in \{(0,0),(0,r),(r(v),r(v)),(r(v),r),(r,r(v)),(r,r)\}.$$
\end{myenu}
\item \bm{ $p'_k\in N(q'_k)\cup N^2(q'_k)\mbox{ and }s'_k\notin N(p'_k)\cup N^2(p'_k)\cup N(q'_k)\cup N^2(q'_k).$}
 \begin{myenu}
 \item If $u=p'_k,$ then
$c_k(u)=(0,r-r(u))$ and $$c_k(v)\in \{(0,r),(r,0),(r(v),r),(r,r(v)),(r,r)\}.$$
\item  If $u=q'_k,$ then $c_k(u)=(0,r)$ and $$c_k(v)\in \{(0,r-r(u)),(r,0),(r(v),r),(r,r(v)),(r,r)\}.$$
\item  If $u=s'_k,$ then $c_k(u)=(r,0)$ and $$c_k(v)\in \{(0,r-r(v)),(0,r),(r(v),r),(r,r(v)),(r,r)\}.$$
\end{myenu}
\item   \bm{$\{p'_k,q'_k,s'_k\}\in \tcal.$}
\begin{myenu}
\item If $u=p'_k,$ then
 $c_k(u)=(-r(u),0)$ and $$c_k(v)\in \{(-r(v),-r(v)),(0,-r(v)),(r(v)-\delta ,r(v)-\delta )\}.$$
\item If $u=q'_k,$ then
 $c_k(u)=(-r(u),-r(u))$ and $$c_k(v)\in \{(-r(v),0),(0,-r(v)),(r(v)-\delta ,r(v)-\delta )\}.$$
\item If $u=s'_k,$ then
$c_k(u)=(0,-r(u))$ and $$c_k(v)\in \{(-r(v),0),(-r(v),-r(v)),(r(v)-\delta ,r(v)-\delta )\}.$$
\end{myenu}
\end{myenu}
In each of these cases,  (\ref{LemEq2}) is satisfied.

\subsection{Proof of (\ref{LemEq3})}

Let $u\in P_k\nd  v\in \displaystyle\bigcup_{l\leq k} P_l \sm \{u\}$ be such that $uv\not\in E(G)$   and $\{u,v\}\subseteq S_w \fos w\in V(G).$
We need to establish that
$$|c_k(v)-c_k(u)|  \geq r(v)+r(u).$$

Since $u\in P_k=\{p_k',q_k',s_k'\},$ which has been chosen at step
\ref{CharVerS43} of   Algorithm \ref{alg41}, it is not an exterior vertex of any star.
So there are no such vertices $u\nd v$ and consequently the above inequality is vacuously true.

%

\subsection{Proof of (\ref{LemEq4})}

Let $u\in P_k$  and  $v\in \displaystyle\bigcup_{l\geq k} P_l\sm \{u\} $ be such that $uv\not\in E(G)$   and $\{u,v\}\not \subseteq S_w \foa w\in V(G).$  We need to establish that
$$|c_k(v)-c_k(u)|  \geq r(v)+r(u).$$
Once $P_k  =\{p_k',q_k',s_k'\}$ is picked and $u\in P_k,$ there exists no $v\in \cup_{l\geq k}P_l$ such that $uv\notin E(G)$
(see algorithm \ref{alg41}). Hence this case is vacuously true.

\subsection{Proof of (\ref{LemEq5})}

Let   $u\in V(G)\nd v\in V(G) \sm \{u\}$ such that $uv\in E(G).$  It is enough to prove that
$$|c_k(v)-c_k(u)|  < r(v)+r(u).$$
We now study this case by case as follows.
\begin{myenu}[label=Case \arabic*.]
\item \bm{$P_k=\{x,y,z\}.$}
Then we have the following cases
\begin{align*}
c_k(u)&\in \{(0,0),(0,r),(r,0),(r(u),r(u)),(r(u),r),(r,r(u)),(r,r)\}\\
\nd  c_k(v)&\in \{(0,0),(0,r),(r,0),(r(v),r(v)),(r(v),r),(r,r(v)),(r,r)\}.
\end{align*}
\item \bm{ $p'_k\in N(q'_k)\cup N^2(q'_k)\mbox{ and }s'_k\notin N(p'_k)\cup N^2(p'_k)\cup N(q'_k)\cup N^2(q'_k).$}
Then
\begin{align*}
c_k(u)&\in \{(0,r-r(u)),(0,r),(r,0),(r(u),r),(r,r(u)),(r,r)\}\\
\nd c_k(v)&\in \{(0,r-r(v)),(0,r),(r,0),(r(v),r),(r,r(v)),(r,r)\}.
\end{align*}
\item   \bm{$\{p'_k,q'_k,s'_k\}\in \tcal.$}
Then we have the following cases
\begin{align*}
c_k(u)&\in \{(-r(u),0),(-r(u),-r(u)),(0,-r(u)),(r(u)-\delta ,r(u)-\delta )\}\\
\nd c_k(v)&\in \{(-r(v),0),(-r(v),-r(v)),(0,-r(v)),(r(v)-\delta ,r(v)-\delta )\}.
\end{align*}
\end{myenu}
In each of these cases,  (\ref{LemEq5}) is satisfied.



\section{\bf Module V: The Case of $\{v_0,w_1,w_2\}$}

Consider the case of independent triplets, picked up through the recursion of step \ref{CharVerS33} of   Algorithm \ref{alg41}. We reserve two dimensions for each such choice. Let $P_k=\{v_0,w_1,w_2\}$ be any such triplet.

Recall that at the recursion of step \ref{CharVerS33} of   Algorithm \ref{alg41},
the collection of unpicked vertices contain vertices from at most one star $\{w\}\cup S_w.$ The center $w$ had already been picked up. We denoted the unpicked vertices of $S_w$ by $D_w.$
Therefore $w_1, w_2\in S_w$ and hence $N(w_1)=\{w\}=N(w_2).$ Consequently, $r(w)=r(w_1)=r(w_2).$


%

Just like  Proposition \ref{propII}, one can establish the following result.
\begin{prop}\label{propV}
If $v\in \ds\bigcup_{i>k}P_i \setminus D_w,$ then $\{vv_0, vw_1, vw_2\}\subset E(G).$
   \end{prop}
Further, for $w_i\in D_w\setminus  \{w_1,w_2\},$ note that $w_iw_1\notin E(G)\nd w_iw_2\notin E(G),$  as these are the exterior vertices of a star. Also note that $r(w_1)=r(w)=r(w_2).$

Based upon our choice of $P_k=\{v_0,w_1,w_2\},$ the corresponding dimensions on the vertices of $G$ are assigned as follows.

\begin{myenu}[label=Step \arabic*.]
\item Define $c_k(w_1):=(-r(w)-\frac{r}{2}+\delta (k+1),r(w)).$

\item Define $c_k(w_2):=(-r(w)-\frac{r}{2}+\delta (k+1),-r(w)).$

\item Define $c_k(v_0):=(r(v_0)+\frac{r}{2}-\delta (k+1),r(v_0)).$

\item Define  $c_k(w):= (-\frac{r}{2}+\delta (k+1),0).$

\item Let $w_j\in N(w)\setminus \{w_1,w_2\}.$ If $w_j\in \bigcup_{i<k}P_i,$ define
$$c_k(w_j):=\bigg(r(w)-\frac{r}{2}+\delta (k+1),-r(w)+\delta\bigg).$$
If $w_j\in \bigcup_{i>k}P_i,$ define
$$c_k(w_j):=\left\{\begin{array}{ll}
(r(w)-\frac{r}{2}+\delta (k+1),-r(w)+\delta) & \mbox{; $w_jv_0\in E(G),$ }\\
(r(w)-\frac{r}{2}+\delta (k+1),-r(w)) & \mbox{; $w_jv_0\notin E(G).$ }
\end{array}\right.$$

\item If $v\in N(v_0)\cup N^2(v_0)\setminus \{v_0\},$ define
 $c_k(v):=(\frac{r}{2}-\delta (k+1),0).$

\item If $v\in V(G)\sm  N_k,$ define
$$c_k(v):=\bigg(r(v)-\frac{r}{2}+\delta (k+1)-\delta ,-r(v)+\delta\bigg).$$
\end{myenu}

\subsection{Proof of (\ref{LemEq1})}

Let $v\in V(G) \nd nv\in N(v).$ We need to establish that
\begin{equation*}
 |c_k(v)-c_k(nv)|  \leq r(v).
\end{equation*}

%
%

Note that for each $u\in V(G),$ the set $N(u)\cup N^2(u)$ is closed under the pseudo-neighborhood operation.
Therefore, it is enough to establish
$$|c_k(nu)-c_k(n^2u)|\leq r(u)\foa  nu\in N(u)\ndfoa n^2u\in N^2(u).$$
Hence we consider the following cases.

\begin{myenu}[label=Case \arabic*.]
\item \bm{$u\in N(v_0).$} Then $m(u)\leq k.$ Note that
$c_k(u)=(\frac{r}{2}-\delta (k+1),0)$   and we have
\\$c_k(nu)\in \{(r(u)+\frac{r}{2}-\delta (k+1),r(u)),(\frac{r}{2}-\delta (k+1),0)\}.$

\item \bm{$u\in N(w_1).$} Then $m(u)\leq k.$  Note that $c_k(u)=(-\frac{r}{2}+\delta (k+1),0)$ and
\\$c_k(nu)\in \{(-r(u)-\frac{r}{2}+\delta (k+1),r(u)),  (-r(u)-\frac{r}{2}+\delta (k+1),-r(u)),\\
(r(u)-\frac{r}{2}+\delta (k+1), -r(u)+\delta ), (r(u)-\frac{r}{2}+\delta (k+1),-r(u))\}.$

\item \bm{$u\in N(w_2).$} Same as the previous  case.

\item \bm{$u\in V(G)\sm  N_k.$} Then  $nu\in V(G)\sm  N_k$ and we have
\\$c_k(u)=(r(u)-\frac{r}{2}+\delta (k+1)-\delta ,-r(u)+\delta )= c_k(nu).$
\end{myenu}

\subsection{Proof of (\ref{LemEq2})}
Let  $u\in P_k$  and  $v\in \displaystyle\bigcup_{l\leq k} P_l\sm \{u\}.$
 We need to  prove that
\begin{equation*}
 |c_k(v)-c_k(u)|  \geq \max\{r(u),r(v)\}.
\end{equation*}

%

Now we consider the following cases.
\begin{myenu}[label=Case \arabic*.]
\item \bm{$u=v_0.$} Then $c_k(u)=(r(v_0)+\frac{r}{2}-\delta (k+1),r(v_0))$ and
\\$c_k(v)\in \{(-r(v)-\frac{r}{2}+\delta (k+1),r(v)),(-r(v)-\frac{r}{2}+\delta (k+1),-r(v)),(-\frac{r}{2}+\delta (k+1),0),\\
(r(v)-\frac{r}{2}+\delta (k+1),-r(v)+\delta ),(\frac{r}{2}-\delta (k+1),0),(r(v)-\frac{r}{2}+\delta (k+1)-\delta ,-r(v)+\delta )\}.$

\item \bm{$u=w_1.$} Then  $c_k(u)=(-r(u)-\frac{r}{2}+\delta (k+1),r(u))$ and
\\$c_k(v)\in \{(-r(v)-\frac{r}{2}+\delta (k+1),-r(v)),(r(v)+\frac{r}{2}-\delta (k+1),r(v)),(-\frac{r}{2}+\delta (k+1),0),\\
(r(v)-\frac{r}{2}+\delta (k+1),-r(v)+\delta ),(\frac{r}{2}-\delta (k+1),0),(r(v)-\frac{r}{2}+\delta (k+1)-\delta ,-r(v)+\delta )\}.$

\item \bm{$u=w_2.$} Then $c_k(u)=(-r(u)-\frac{r}{2}+\delta (k+1),-r(u))$ and
\\$c_k(v)\in \{(-r(v)-\frac{r}{2}+\delta (k+1),r(v)),(r(v)+\frac{r}{2}-\delta (k+1),r(v)),(-\frac{r}{2}+\delta (k+1),0),\\
(r(v)-\frac{r}{2}+\delta (k+1),-r(v)+\delta ),(\frac{r}{2}-\delta (k+1),0),(r(v)-\frac{r}{2}+\delta (k+1)-\delta ,-r(v)+\delta )\}.$

\end{myenu}

%

\subsection{Proof of (\ref{LemEq3})}

Let  $u\in P_k\nd  v\in  \bigcup_{l\leq k} P_l \sm \{u\}$  such that $uv\not\in E(G)$   and $\{u,v\}\subseteq S_w \fos w\in V(G).$\\
 We need to  prove that
$$|c_k(v)-c_k(u)|  \geq r(v)+r(u).$$

We have the following cases.

\begin{myenu}[label=Case \arabic*.]
\item \bm{$u=v_0.$} This case is impossible. Because when $P_k$ was picked up, at that time the set of unpicked vertices had elements at most from one star $S_w.$ We denoted that subset of $S_w$ by $D_w$ and we had chosen  $v_0$ outside $ D_w.$

\item \bm{$u=w_1.$} Then $c_k(u)=(-r(u)-\frac{r}{2}+\delta (k+1),r(u))$ and
\\$c_k(v)\in \{(-r(v)-\frac{r}{2}+\delta (k+1),-r(v)),(r(v)-\frac{r}{2}+\delta (k+1),-r(v)+\delta )\}.$

\item \bm{$u=w_2.$} Then $c_k(u)=(-r(u)-\frac{r}{2}+\delta (k+1),-r(u))$ and
\\$c_k(v)\in \{(-r(v)-\frac{r}{2}+\delta (k+1),r(v)),(r(v)-\frac{r}{2}+\delta (k+1),-r(v)+\delta )\}.$
\end{myenu}

\subsection{Proof of (\ref{LemEq4})}

Let  $u\in P_k$  and  $v\in \displaystyle\bigcup_{l\geq k} P_l\sm \{u\} $ be such that $uv\not\in E(G)$   and $\{u,v\}\not \subseteq S_w \foa w\in V(G).$\\
We need to  prove that
$$|c_k(v)-c_k(u)|  \geq r(v)+r(u).$$

Now, we consider the following cases.
\begin{myenu}[label=Case \arabic*.]
\item \bm{$u=v_0.$} Then $c_k(u)=(r(u)+\frac{r}{2}-\delta (k+1),r(u))$ and
\\$c_k(v)\in \{(-r(v)-\frac{r}{2}+\delta (k+1),r(v)),(-r(v)-\frac{r}{2}+\delta (k+1),-r(v)),(r(v)-\frac{r}{2}+\delta (k+1),-r(v))\}.$

\item \bm{$u=w_1.$} Then $c_k(u)=(-r(u)-\frac{r}{2}+\delta (k+1),r(u))\nd c_k(v)=(r(v)+\frac{r}{2}-\delta (k+1),r(v)).$

\item \bm{$u=w_2.$} Then $c_k(u)=(-r(u)-\frac{r}{2}+\delta (k+1),-r(u))\nd c_k(v)=(r(v)+\frac{r}{2}-\delta (k+1),r(v)).$

\end{myenu}

 \subsection{Proof of (\ref{LemEq5})}

Let   $u\in V(G)\nd v\in V(G) \sm \{u\}$ such that $uv\in E(G).$  It is enough to prove that
$$|c_k(v)-c_k(u)|  < r(v)+r(u).$$

%
%
Recall that $P_k=\{v_0,w_1,w_2\}.$
Suppose that $u \in P_i\nd v \in P_j.$ Then
  $$r(u) \geq r-2\delta i \mbox{ and } r(v) \geq r-2\delta j.$$
\begin{myenu}[label=Case \arabic*.]
\item \bm{$k\notin \{i,j\}.$} Let $Q:=(0,-\delta )\nd A(y)$ denote the set
\\$\{(-\frac{r}{2}+\delta (k+1),0),(r(y)-\frac{r}{2}+\delta (k+1),-r(y)+\delta ),(r(y)-\frac{r}{2}+\delta (k+1),-r(y)+\delta ),\\
(r(y)-\frac{r}{2}+\delta (k+1),-r(y)),(\frac{r}{2}-\delta (k+1),0),(r(y)-\frac{r}{2}+\delta (k+1)-\delta ,-r(y)+\delta )\}.$

Then $c_k(u)\in A(u)\nd c_k(v)\in A(v).$ Note that we have
$$|c_k(u)-Q|<r(u)\nd |c_k(v)-Q|<r(v).$$
Hence $|c_k(u)-c_k(v)|\leq |c_k(u)-Q|+|Q-c_k(v)|<r(u)+r(v).$

\item \bm{$k\in \{i,j\}.$} Without loss of generality, we assume that $k=i.$
\begin{myenu}
\item \bm{$u=v_0.$} Then $c_k(u)=(r(u)+\frac{r}{2}-\delta (k+1),r(u))$ and
\\$c_k(v)\in \{(-\frac{r}{2}+\delta (k+1),0),(r(v)-\frac{r}{2}+\delta (k+1),-r(v)+\delta ),(r(v)-\frac{r}{2}+\delta (k+1),\\
-r(v)+\delta ),(\frac{r}{2}-\delta (k+1),0),(r(v)-\frac{r}{2}+\delta (k+1)-\delta ,-r(v)+\delta )\}.$

\item \bm{$u=w_1.$} Then $c_k(u)=(-r(u)-\frac{r}{2}+\delta (k+1),r(u))$ and
\\$c_k(v)\in \{(-\frac{r}{2}+\delta (k+1),0),(\frac{r}{2}-\delta (k+1),0),(r(v)-\frac{r}{2}+\delta (k+1)-\delta ,-r(v)+\delta )\}.$

\item \bm{$u=w_2.$} Then $c_k(u)=(-r(u)-\frac{r}{2}+\delta (k+1),-r(u))$ and
\\$c_k(v)\in \{(-\frac{r}{2}+\delta (k+1),0),(\frac{r}{2}-\delta (k+1),0),(r(v)-\frac{r}{2}+\delta (k+1)-\delta ,-r(v)+\delta )\}.$
\end{myenu}
\end{myenu}

\section{\bf Module VI: The Case of $\{w_0, v_1, v_2\}$}

Consider the case of triplets, picked up through the recursion of   step \ref{CharVerS35} of   Algorithm \ref{alg41}. We reserve two dimensions for each such choice. Let $P_k=\{w_0,v_1,v_2\}$ be any such triplet,
 with only $w_0\in D_w, v_1v_2\in E(G)$ and exactly two edges among them.
Then exactly one of $w_0v_1\in E(G)$ and  $w_0v_2 \in E(G) $ is true. Without loss of generality, assume that
$w_0v_1\in E(G)\nd w_0v_2\notin E(G).$
We further divide it into two  subcases given by
 $$v_2\notin N(v_1)\cup N^2(v_1)\nd v_2\in N(v_1)\cup N^2(v_1).$$
We will treat these cases separately, while assigning   $c_k$ on the vertices of $G.$
\subsection{\bf Subcase VI a}

If $v_2\notin N(v_1)\cup N^2(v_1),$
the corresponding dimensions  are assigned as follows.


\begin{myenu}[label=Step \arabic*.]
\item Define $c_k(w_0):=(r(w_0)-\delta ,r(w_0)+\frac{r}{2}-\delta (k+1)).$

\item Define $c_k(v_1):=(-r(v_1),r(v_1)-\frac{r}{2}+\delta (k+1)-\delta ).$

\item Define $c_k(v_2):=(r(v_2)-\delta ,-r(v_2)-\frac{r}{2}+\delta (k+1)).$

\item Let $v\in V(G)\sm  N_k.$ If $v\in \cup_{k<i} P_k,$ define $c_k(v):=(r(v)-\delta ,0).$
\\If $v\in \cup_{k>i} P_k,$ define
$$c_k(v):=\left\{\begin{array}{ll}
(r(v)-\delta ,0) & \mbox{; $vv_1\in E(G),$ }\\
(r(v),0) & \mbox{; $vv_1\notin E(G).$ }
\end{array}\right.$$

\item If $v\in N(v_1)\cup N^2(v_1)\sm \{v_1\},$ define $c_k(v):=(0,0).$

\item Let $v\in N(v_2)\cup N^2(v_2)\sm \{v_2\}.$
If $v\in \cup_{i<k}P_i,$
define $c_k(v):=(r(v)-\delta ,-\frac{r}{2}+\delta (k+1)).$

If $v\in \cup_{i>k}P_i,$
define
$$c_k(v):=\left\{\begin{array}{ll}
(r(v)-\delta ,-\frac{r}{2}+\delta (k+1)) & \mbox{; $vv_1\in E(G),$ }\\
(r(v),-\frac{r}{2}+\delta (k+1)) & \mbox{; $vv_1\notin E(G).$ }
\end{array}\right.$$

\item Let $v\in N(w_0).$ Note that $v\in \cup_{i<k}P_i.$ Define
$$c_k(v):=(r-2\delta (k+1),\frac{r}{2}-\delta (k+1))$$

\item Let $v\in N^2(w_0).$  If $v\in \cup_{i<k}P_i,$ define
$$c_k(v):=(-r(v)+r-2\delta (k+1),-r(v)+\frac{r}{2}-\delta (k+1)).$$
If $v\in \cup_{i>k}P_i,$
define $$c_k(v):=\left\{\begin{array}{ll}
(r(v)-\delta ,\frac{r}{2}-\delta (k+1)) & \mbox{; $vv_1\in E(G),$ }\\
(r(v),\frac{r}{2}-\delta (k+1)) & \mbox{; $vv_1\notin E(G).$ }
\end{array}\right.$$

\end{myenu}

\subsubsection{Proof of (\ref{LemEq1})}

Let $v\in V(G) \nd nv\in N(v).$ We need to establish that
\begin{equation*}
 |c_k(v)-c_k(nv)|  \leq r(v).
\end{equation*}

Note that for each $u\in V(G),$ the set $N(u)\cup N^2(u)$ is closed under the pseudo-neighborhood operation.
Therefore, it is enough to establish
$$|c_k(nu)-c_k(n^2u)|\leq r(u)\foa  nu\in N(u)\ndfoa n^2u\in N^2(u).$$
Also, note that if $u\in N_k,$ then $m(u)\leq k$ and thus $r(u)\geq r-2\del k.$
Hence we consider the following cases.



\begin{myenu}[label=Case \arabic*.]
\item \bm{$u\in N(w_0).$} Then   $c_k(u)=(r-2\delta (k+1),\frac{r}{2}-\delta (k+1))$ we have
\\$c_k(nu)\in \{(r(u)-\delta ,r(u)+\frac{r}{2}-\delta (k+1)),(-r(u)+r-2\delta (k+1),-r(u)+\frac{r}{2}-\delta (k+1)),\\
(r(u)-\delta ,\frac{r}{2}-\delta (k+1)),(r(u),\frac{r}{2}-\delta (k+1))\}.$

\item \bm{$u\in N(v_1).$} Then    $c_k(u)=(0,0)\nd c_k(nu)\in \{(-r(u),r(u)-\frac{r}{2}+\delta (k+1)-\delta ),(0,0)\}.$

\item \bm{$u\in N(v_2).$} Then   $c_k(u)\in \{(r(u)-\delta ,-\frac{r}{2}+\delta (k+1)),(r(u),-\frac{r}{2}+\delta (k+1))\}$ and \\$c_k(nu)\in \{(r(u)-\delta ,-r(u)-\frac{r}{2}+\delta (k+1)),(r(u)-\delta ,-\frac{r}{2}+\delta (k+1)),(r(u),-\frac{r}{2}+\delta (k+1))\}.$

\item \bm{$u\in V(G)\sm  N_k.$} Then  $nu\in V(G)\sm  N_k$ and both $c_k(u), c_k(nu)\in \{(r(u)-\delta ,0),(r(u),0)\}.$
\end{myenu}

\subsubsection{Proof of (\ref{LemEq2})}

Let  $u\in P_k$  and  $v\in \displaystyle\bigcup_{l\leq k} P_l\sm \{u\}.$
 We need to  prove that
\begin{equation*}
 |c_k(v)-c_k(u)|  \geq \max\{r(u),r(v)\}.
\end{equation*}


Now we consider the following cases.
\begin{myenu}[label=Case \arabic*.]
\item \bm{$u=w_0.$} Then $c_k(u)=(r(u)-\delta ,r(u)+\frac{r}{2}-\delta (k+1))$ and
\\$c_k(v)\in \{(-r(v),r(v)-\frac{r}{2}+\delta (k+1)-\delta ),(r(v)-\delta ,-r(v)-\frac{r}{2}+\delta (k+1)),(r(v)-\delta ,0),(0,0),(r(v)-\delta ,-\frac{r}{2}+\delta (k+1)),(r-2\delta (k+1),\frac{r}{2}-\delta (k+1)),(-r(v)+r-2\delta (k+1),-r(v)+\frac{r}{2}-\delta (k+1))\}.$

\item \bm{$u=v_1.$} Then $c_k(u)=(-r(u),r(u)-\frac{r}{2}+\delta (k+1)-\delta )$ and
\\$c_k(v)\in \{(r(v)-\delta ,r(v)+\frac{r}{2}-\delta (k+1)),(r(v)-\delta ,-r(v)-\frac{r}{2}+\delta (k+1)),(r(v)-\delta ,0),(0,0),(r(v)-\delta ,-\frac{r}{2}+\delta (k+1)),(r-2\delta (k+1),\frac{r}{2}-\delta (k+1)),(-r(v)+r-2\delta (k+1),-r(v)+\frac{r}{2}-\delta (k+1))\}.$

\item \bm{$u=v_2.$} Then $c_k(u)=(r(u)-\delta ,-r(u)-\frac{r}{2}+\delta (k+1))$ and
\\$c_k(v)\in \{(r(v)-\delta ,r(v)+\frac{r}{2}-\delta (k+1)),(-r(v),r(v)-\frac{r}{2}+\delta (k+1)-\delta ),(r(v)-\delta ,0),(0,0),(r(v)-\delta ,-\frac{r}{2}+\delta (k+1)),(r-2\delta (k+1),\frac{r}{2}-\delta (k+1)),(-r(v)+r-2\delta (k+1),-r(v)+\frac{r}{2}-\delta (k+1))\}.$
\end{myenu}

\subsubsection{Proof of (\ref{LemEq3})}

Let  $u\in P_k\nd  v\in \displaystyle\bigcup_{l\leq k} P_l \sm \{u\}$  such that $uv\not\in E(G)$   and $\{u,v\}\subseteq S_w \fos w\in V(G).$
 We need to  prove that
$$|c_k(v)-c_k(u)|  \geq r(v)+r(u).$$


The following are the only cases.

\begin{myenu}[label=Case \arabic*.]
\item \bm{$u=w_0.$} Then $c_k(u)=(r(u)-\delta ,r(u)+\frac{r}{2}-\delta (k+1))$ and
\\$c_k(v)=(-r(v)+r-2\delta (k+1),-r(v)+\frac{r}{2}-\delta (k+1)).$
\item \bm{$u=v_1.$}  This case fails.
\item \bm{$u=v_2.$} This case fails.
\end{myenu}

\subsubsection{Proof of (\ref{LemEq4})}

Let  $u\in P_k$  and  $v\in \displaystyle\bigcup_{l\geq k} P_l\sm \{u\} $ be such that $uv\not\in E(G)$   and $\{u,v\}\not \subseteq S_w \foa w\in V(G).$
\\We need to  prove that
$$|c_k(v)-c_k(u)|  \geq r(v)+r(u).$$

%
%

\begin{myenu}[label=Case \arabic*.]
\item \bm{$u=w_0.$} Then either $uv\in E(G)\mbox{ or }\{u,v\}\subseteq S_w$ for some $w\in V(G)\mbox{ or }u=w_0\nd v=v_2.$
\\Hence $u=w_0\nd v=v_2.$ \Csq $c_k(u)=(r(u)-\delta ,r(u)+\frac{r}{2}-\delta (k+1))$ and
\\$c_k(v)=(r(v)-\delta ,-r(v)-\frac{r}{2}+\delta (k+1)).$

\item \bm{$u=v_2.$} Then either $uv\in E(G)\mbox{ or }\{u,v\}\subseteq S_w$ for some $w\in V(G)\mbox{ or }u=v_2\nd v=w_0.$
\\Therefore $u=v_2\nd v=w_0.$ Hence $c_k(u)=(r(u)-\delta ,-r(u)-\frac{r}{2}+\delta (k+1))$ and
\\$c_k(v)=(r(v)-\delta ,r(v)+\frac{r}{2}-\delta (k+1)).$

\item \bm{$u=v_1.$} Then $c_k(u)=(-r(u),r(u)-\frac{r}{2}+\delta (k+1)-\delta )$ and
\\$c_k(v)\in \{(r(v),0),(r(v),-\frac{r}{2}+\delta (k+1)),(r(v),\frac{r}{2}-\delta (k+1))\}.$
\end{myenu}

\subsubsection{Proof of (\ref{LemEq5})}

Let   $u\in V(G)\nd v\in V(G) \sm \{u\}$ such that $uv\in E(G).$  It is enough to prove that
$$|c_k(v)-c_k(u)|  < r(v)+r(u).$$

%
%

Assume that  one of $u\nd v$ belongs to $P_i$ and other belongs to $P_j.$
We have the following cases.

\begin{myenu}[label=Case \arabic*.]
\item \bm{$k\notin \{i,j\}.$}  Let $Q:=(\delta ,0)\nd A(a)$ denote the set
\\$\{(r(a)-\delta ,0),(r(a),0),(0,0),(r(a)-\delta ,-\frac{r}{2}+\delta (k+1)),(r(a),-\frac{r}{2}+\delta (k+1)),(r-2\delta (k+1),\frac{r}{2}-\delta (k+1)),(-r(a)+r-2\delta (k+1),-r(a)+\frac{r}{2}-\delta (k+1)),(r(a)-\delta ,\frac{r}{2}-\delta (k+1)),(r(a),\frac{r}{2}-\delta (k+1))\}.$
\\Then $c_k(u)\in A(u)\nd c_k(v)\in A(v).$
Note that $|c_k(u)-Q|<r(u)\nd |c_k(v)-Q|<r(v).$ Hence $|c_k(u)-c_k(v)|<r(u)+r(v).$

\item \bm{$k\in \{i,j\}.$} Without loss of generality, let $k=i.$
\begin{myenu}
\item \bm{$u=w_0.$} Then $c_k(u)=(r(u)-\delta ,r(u)+\frac{r}{2}-\delta (k+1))$ and
\\$c_k(v)\in \{(-r(v),r(v)-\frac{r}{2}+\delta (k+1)-\delta ),(r(v)-\delta ,0),(r(v),0),(0,0),(r(v)-\delta ,-\frac{r}{2}+\delta (k+1)),(r(v),-\frac{r}{2}+\delta (k+1)),(r-2\delta (k+1),\frac{r}{2}-\delta (k+1))\}.$
\item \bm{$u=v_1.$} Then $c_k(u)=(-r(u),r(u)-\frac{r}{2}+\delta (k+1)-\delta )$ and
\\$c_k(v)\in \{(r(v)-\delta ,r(v)+\frac{r}{2}-\delta (k+1)),(r(v)-\delta ,-r(v)-\frac{r}{2}+\delta (k+1)),(r(v)-\delta ,0),(0,0),(r(v)-\delta ,-\frac{r}{2}+\delta (k+1)),(r-2\delta (k+1),\frac{r}{2}-\delta (k+1)),(-r(v)+r-2\delta (k+1),-r(v)+\frac{r}{2}-\delta (k+1)),(r(v)-\delta ,\frac{r}{2}-\delta (k+1))\}.$
\item \bm{$u=v_2.$} Then $c_k(u)=(r(u)-\delta ,-r(u)-\frac{r}{2}+\delta (k+1))$ and
\\$c_k(v)\in \{(-r(v),r(v)-\frac{r}{2}+\delta (k+1)-\delta ),(r(v)-\delta ,0),(r(v),0),(0,0),(r(v)-\delta ,-\frac{r}{2}+\delta (k+1)),(r(v),-\frac{r}{2}+\delta (k+1)),(r-2\delta (k+1),\frac{r}{2}-\delta (k+1)),(-r(v)+r-2\delta (k+1),-r(v)+\frac{r}{2}-\delta (k+1)),(r(v)-\delta ,\frac{r}{2}-\delta (k+1)),(r(v),\frac{r}{2}-\delta (k+1))\}.$
\end{myenu}
\end{myenu}

\subsection{\bf Subcase VI b}

Now, we consider the subcase of $P_k=\{w_0, v_1, v_2\},$ with $w_0v_1, v_1v_2\in E(G)$ when  $v_2\in N(v_1)\cup N^2(v_1).$
Recall that  $w_0\in D_w$ and there are exactly two edges among  $w_0, v_1 \nd  v_2..$
The corresponding dimensions are assigned as follows.

\begin{myenu}[label=Case \arabic*.]\item Define $c_k(w_0):=(r(w_0)-\delta ,r(w_0)+\frac{r}{2}-\delta (k+1)).$

\item Define $c_k(v_1):=(-r(v_1),-\frac{r}{2}+\delta (k+1)).$

\item Define $c_k(v_2):=(-r(v_2)+\delta ,-r(v_2)-\frac{r}{2}+\delta (k+1)).$

\item Let $v\in V(G)\sm  N_k.$ If $v\in \cup_{k<i} P_k,$
define $c_k(v):=(r(v)-\delta ,0).$

Suppose  $v\in \cup_{k>i} P_k.$ Note that in this case $vv_2\notin E(G).$ Define $c_k(v):=(r(v),0).$

\item Let $v\in N(v_1)\cup N^2(v_1)\cup N(v_2)\cup N^2(v_2)\setminus \{v_1,v_2\}.$ Note that $v\in T,$ for some
 $T\in \tcal.$ Define $$c_k(v):=\bigg(0,-\frac{r}{2}+\delta (k+1)\bigg).$$

\item Let $v\in N(w_0).$ Note that $v\in \cup_{i<k}P_i.$ Define
 $c_k(v):=(r(v)-\delta ,\frac{r}{2}-\delta (k+1))$

\item Let $v\in N^2(w_0)\sm \{w_0\}.$
If $v\in \cup_{i<k}P_i,$ define
$$c_k(v):=\bigg(r(v)-\delta ,-r(v)+\frac{r}{2}-\delta (k+1)\bigg).$$
If $v\in \cup_{i>k}P_i,$ define
$$c_k(v):=\left\{\begin{array}{ll}
(r(v)-\delta ,\frac{r}{2}-\delta (k+1)) & \mbox{; $vv_2\in E(G),$ }\\
(r(v),\frac{r}{2}-\delta (k+1)) & \mbox{; $vv_2\notin E(G).$ }
\end{array}\right.$$
\end{myenu}


\subsubsection{Proof of (\ref{LemEq1})}

Let $v\in V(G) \nd nv\in N(v).$ We need to establish that
\begin{equation*}
 |c_k(v)-c_k(nv)|  \leq r(v).
\end{equation*}

Note that for each $u\in V(G),$ the set $N(u)\cup N^2(u)$ is closed under the pseudo-neighborhood operation.
Therefore, it is enough to establish
$$|c_k(nu)-c_k(n^2u)|\leq r(u)\foa  nu\in N(u)\ndfoa n^2u\in N^2(u).$$
Hence we consider the following cases.

%


\begin{myenu}[label=Case \arabic*.]
\item \bm{$u\in N(w_0).$} Then
$c_k(u)=(r(u)-\delta ,\frac{r}{2}-\delta (k+1))$ and
\\$c_k(nu)\in \{(r(u)-\delta ,r(u)+\frac{r}{2}-\delta (k+1)),(r(u)-\delta ,-r(u)+\frac{r}{2}-\delta (k+1)),(r(u)-\delta ,\frac{r}{2}-\delta (k+1)),(r(u),\frac{r}{2}-\delta (k+1))\}.$

\item \bm{$u\in N(v_1).$}
  This case also covers the case $u\in N(v_2).$ Because $v_2\in N(v_1)\cup N^2(v_1).$ Now
\\$c_k(u)\in \{(-r(u)+\delta ,-r(u)-\frac{r}{2}+\delta (k+1)),(0,-\frac{r}{2}+\delta (k+1))\}$ and
\\$c_k(nu)\in \{(-r(u),-\frac{r}{2}+\delta (k+1)),(-r(u)+\delta ,-r(u)-\frac{r}{2}+\delta (k+1)),(0,-\frac{r}{2}+\delta (k+1))\}.$

\item \bm{$u\in V(G)\sm  N_k.$} Then  $nu\in V(G)\sm  N_k$ and both
\\$c_k(u),c_k(nu)\in \{(r(u)-\delta ,0),(r(u),0)\}.$
\end{myenu}

\subsubsection{Proof of (\ref{LemEq2})}

Let  $u\in P_k$  and  $v\in \displaystyle\bigcup_{l\leq k} P_l\sm \{u\}.$
 We need to  prove that
\begin{equation*}
 |c_k(v)-c_k(u)|  \geq \max\{r(u),r(v)\}.
\end{equation*}

%

We have the following cases.
\begin{myenu}[label=Case \arabic*.]
\item \bm{$u=w_0.$} Then $c_k(u)=(r(u)-\delta ,r(u)+\frac{r}{2}-\delta (k+1))$ and
\\$c_k(v)\in \{(-r(v),-\frac{r}{2}+\delta (k+1)),(-r(v)+\delta ,-r(v)-\frac{r}{2}+\delta (k+1)),(r(v)-\delta ,0),\\(0,-\frac{r}{2}+\delta (k+1)),(r(v)-\delta ,\frac{r}{2}-\delta (k+1)),(r(v)-\delta ,-r(v)+\frac{r}{2}-\delta (k+1))\}.$

\item \bm{$u=v_1.$} Then $c_k(u)=(-r(u),-\frac{r}{2}+\delta (k+1))$ and
\\$c_k(v)\in \{(r(v)-\delta ,r(v)+\frac{r}{2}-\delta (k+1)),(-r(v)+\delta ,-r(v)-\frac{r}{2}+\delta (k+1)),(r(v)-\delta ,0),\\
(0,-\frac{r}{2}+\delta (k+1)),(r(v)-\delta ,\frac{r}{2}-\delta (k+1)),(r(v)-\delta ,-r(v)+\frac{r}{2}-\delta (k+1))\}.$

\item \bm{$u=v_2.$} Then $c_k(u)=(-r(u)+\delta ,-r(u)-\frac{r}{2}+\delta (k+1))$ and
\\$c_k(v)\in \{(r(v)-\delta ,r(v)+\frac{r}{2}-\delta (k+1)),(-r(v),-\frac{r}{2}+\delta (k+1)),(r(v)-\delta ,0),\\
(0,-\frac{r}{2}+\delta (k+1)),(r(v)-\delta ,\frac{r}{2}-\delta (k+1)),(r(v)-\delta ,-r(v)+\frac{r}{2}-\delta (k+1))\}.$
\end{myenu}

\subsubsection{Proof of (\ref{LemEq3})}

Let  $u\in P_k\nd  v\in \displaystyle\bigcup_{l\leq k} P_l \sm \{u\}$  such that $uv\not\in E(G)$   and $\{u,v\}\subseteq S_w \fos w\in V(G).$
 \\ We need to  prove that
$$|c_k(v)-c_k(u)|  \geq r(v)+r(u).$$

%
%
We have the following cases.

\begin{myenu}[label=Case \arabic*.]\item \bm{$u=w_0.$} Then
\\$c_k(u)=(r(u)-\delta ,r(u)+\frac{r}{2}-\delta (k+1))\nd c_k(v)=(r(v)-\delta ,-r(v)+\frac{r}{2}-\delta (k+1)).$
\item \bm{$v=v_1.$}  This case is impossible.
\item \bm{$v=v_2.$} This case impossible.
\end{myenu}

\subsubsection{Proof of (\ref{LemEq4})}

Let  $u\in P_k$  and  $v\in \displaystyle\bigcup_{l\geq k} P_l\sm \{u\} $ be such that $uv\not\in E(G)$   and $\{u,v\}\not \subseteq S_w \foa w\in V(G).$
\\ We need to  prove that
$$|c_k(v)-c_k(u)|  \geq r(v)+r(u).$$

%
%

We have the following cases.
\begin{myenu}[label=Case \arabic*.]\item \bm{$u=w_0.$} Then either $uv\in E(G)\mbox{ or }\{u,v\}\subseteq S_w$ for some $w\in V(G)\mbox{ or }u=w_0\nd v=v_2.$
\\Hence $u=w_0\nd v=v_2.$ Therefore
\\$c_k(u)=(r(u)-\delta ,r(u)+\frac{r}{2}-\delta (k+1))\nd c_k(v)=(-r(v)+\delta ,-r(v)-\frac{r}{2}+\delta (k+1)).$

\item \bm{$u=v_2.$} Then either $uv\in E(G)\mbox{ or }\{u,v\}\subseteq S_w$ for some $w\in V(G)\mbox{ or }u=v_2\nd v=w_0.$
\\Hence $u=v_2\nd v=w_0.$ Therefore
\\$c_k(u)=(-r(u)+\delta ,-r(u)-\frac{r}{2}+\delta (k+1))\nd c_k(v)=(r(v)-\delta ,r(v)+\frac{r}{2}-\delta (k+1)).$

\item \bm{$u=v_1.$} Then $c_k(u)=(-r(u),-\frac{r}{2}+\delta (k+1))$  and $c_k(v)\in \{(r(v),0),(r(v),\frac{r}{2}-\delta (k+1))\}.$
\end{myenu}

\subsubsection{Proof of (\ref{LemEq5})}

Let   $u\in V(G)\nd v\in V(G) \sm \{u\}$ such that $uv\in E(G).$  It is enough to prove that
$$|c_k(v)-c_k(u)|  < r(v)+r(u).$$

%
%

Assume that $u\in P_i\nd v\in P_j.$
  \begin{myenu}[label=Case \arabic*.]
  \item \bm{$k\notin \{i,j\}.$} Let $Q:=(\delta, 0)\nd A(a)$ denote the set
\\$\{(r(a)-\delta ,0),(r(a),0),(0,-\frac{r}{2}+\delta (k+1)),(r(a)-\delta ,\frac{r}{2}-\delta (k+1)),\\
(r(a)-\delta ,-r(a)+\frac{r}{2}-\delta (k+1)),(r(a)-\delta ,\frac{r}{2}-\delta (k+1)),(r(a),\frac{r}{2}-\delta (k+1))\}.$
%
\\ Then $c_k(u)\in A(u)\nd c_k(v)\in A(v).$ Further note that
$$|c_k(u)-Q|<r(u)\nd |c_k(v)-Q|<r(v).$$
Hence $|c_k(u)-c_k(v)|\leq |c_k(u)-Q|+|Q-c_k(v)|<r(u)+r(v).$

\item \bm{$k\in \{i,j\}.$} Without loss of generality, assume that $k=i.$
\begin{myenu}
\item \bm{$u=w_0.$} Then $c_k(u)=(r(u)-\delta ,r(u)+\frac{r}{2}-\delta (k+1))\nd c_k(v)$ belongs to
\\$  \{(-r(v),-\frac{r}{2}+\delta (k+1)),(r(v)-\delta ,0),(r(v),0),(0,-\frac{r}{2}+\delta (k+1)),(r(v)-\delta ,\frac{r}{2}-\delta (k+1))\}.$
\item \bm{$u=v_1.$} Then $c_k(u)=(-r(u),-\frac{r}{2}+\delta (k+1))\nd c_k(v)$ belongs to
\\$ \{(r(v)-\delta ,r(v)+\frac{r}{2}-\delta (k+1)),(-r(v)+\delta ,-r(v)-\frac{r}{2}+\delta (k+1)),(r(v)-\delta ,0),(0,-\frac{r}{2}+\delta (k+1)),(r(v)-\delta ,\frac{r}{2}-\delta (k+1)),(r(v)-\delta ,-r(v)+\frac{r}{2}-\delta (k+1)),(r(v)-\delta ,\frac{r}{2}-\delta (k+1))\}.$
\item \bm{$u=v_2.$} Then $c_k(u)=(-r(u)+\delta ,-r(u)-\frac{r}{2}+\delta (k+1))\nd c_k(v)$ belongs to
\\$ \{(-r(v),-\frac{r}{2}+\delta (k+1)),(r(v)-\delta ,0),(r(v),0),(0,-\frac{r}{2}+\delta (k+1)),
(r(v)-\delta ,\frac{r}{2}-\delta (k+1)),(r(v)-\delta ,-r(v)+\frac{r}{2}-\delta (k+1)),(r(v)-\delta,
\frac{r}{2}-\delta (k+1)),(r(v),\frac{r}{2}-\delta (k+1))\}.$ \qedhere
\end{myenu}
\end{myenu}

%
%
%
%
%

\section{\bf Module VII: Triplets having exactly one edge}

Consider the triplets of vertices picked up through the recursion of    step   \ref{CharVerS24} of Algorithm \ref{alg41}.
We reserve two dimensions for each such triplet.
Let $P_k:=\{p_k,q_k,s_k\}$  such  a triplet with   exactly one edge among its vertices.
Without loss of generality, assume that $p_kq_k\in E(G).$
Below, we assign  the corresponding  dimensions $c_k$ on the vertices of $G.$

\begin{myenu} [label=Case \arabic*.]
\item \textbf{Either }\bm{$p_k\in N(q_k)\mbox{ or }p_k\in N^2(q_k).$} Then we note
that  $p_k$ and $q_k$ do not belong to any star of our factor.
 Therefore, either $p_kq_k$ is a matching edge (that is $pq\in M$) or $\{p_k,q_k,t_k\}\in \tcal \fos t_k\in V(G).$
\begin{myenu} [label=Step \arabic*.]
\item Define $c_k(p_k):=(-r+2\delta(k+1)-r(p_k),-r+2\delta(k+1)).$

\item Define $c_k(q_k):=(-r+2\delta(k+1)-\delta ,-r+2\delta(k+1)-r(q_k)).$

\item Define $c_k(s_k):=(r(s_k),r(s_k)).$

\item If $\{p_k,q_k,v\}\in \tcal\fos v\in V(G),$ recall that we denote this $v$ by $t_k.$ Define $c_k(t_k)$ as follows:
\begin{myenu} [label=Case \arabic*.]
 \item If $t_k\in \cup_{j<k}P_j,$ define $c_k(t_k):=(-r+2\delta(k+1),-r+2\delta(k+1)).$
\item If $t_k\notin \cup_{j<k}P_j,$ we define
$$c_k(t_k):=\left\{\begin{array}{ll}
(-r+2\delta(k+1),-r+2\delta(k+1)) & \mbox{; $t_ks_k\in E(G),$  }\\
(-r(t_k),-r(t_k))		  & \mbox{; $t_ks_k\notin E(G).$  }
\end{array}\right.$$
\end{myenu}
\item Let $ns_k\in N(s_k).$
\begin{myenu}[label=Case \arabic*.]
\item If $ns_k\in \cup_{j<k}P_j,$ define $c_k(ns_k):=(0,0).$
\item If $ns_k\notin \cup_{j<k}P_j,$ we define $c_k(ns_k)$ as follows
$$\left\{\begin{array}{ll}
(0,0) & \mbox{; $ns_kp_k\in E(G), ns_kq_k\in E(G),$}\\
(-r+2\delta(k+1)+r(s_k),0) & \mbox{; $ns_kp_k\notin E(G), ns_kq_k\in E(G),$}\\
(0,-r+2\delta(k+1)+r(s_k)) & \mbox{; $ns_kp_k\in E(G), ns_kq_k\notin E(G),$}\\
(-r+2\delta(k+1)+r(s_k),\\
 \hspace{.3in} -r+2\delta(k+1)+r(s_k)) & \mbox{; $ns_kp_k\notin E(G), ns_kq_k\notin E(G).$}
\end{array}\right.$$
 \end{myenu}
\item Let $n^2s_k\in N(ns_k)\setminus \{s_k\}.$
\begin{myenu}[label=Case \arabic*.]
\item If $n^2s_k\in \cup_{j<k}P_j$ and  $n^2s_k\in S_u\fos u,$ then define
$$c_k(n^2s_k):=(-r+2\delta(k+1)+r(s_k)-\delta ,-r(s_k)).$$
\item If $n^2s_k\in \cup_{j<k}P_j$  and $n^2s_k\notin S_u\foa$stars $S_u,$ then define
$$c_k(n^2s_k):=(0,0).$$
\item If $n^2s_k\notin \cup_{j<k}P_j,$ then define   $c_k(n^2s_k)$ as follows:
$$ \left\{\begin{array}{ll}
(0,0) & \mbox{; $n^2s_kp_k\in E(G), n^2s_kq_k\in E(G),$}\\
(-r+2\delta(k+1)+r(s_k),0) & \mbox{; $n^2s_kp_k\notin E(G), n^2s_kq_k\in E(G),$}\\
(0,-r+2\delta(k+1)+r(s_k)) & \mbox{; $n^2s_kp_k\in E(G), n^2s_kq_k\notin E(G),$}\\
(-r+2\delta(k+1)+r(s_k),\\
 \hspace{.3in} -r+2\delta(k+1)+r(s_k)) & \mbox{; $n^2s_kp_k\notin E(G), n^2s_kq_k\notin E(G).$}
\end{array}\right.$$
 \end{myenu}
\item Let $v\in V(G)\setminus V_k.$ Consider the following cases:
If $v\in \cup_{j<k}P_j$ then define
$$c_k(v):=(-r/2+\delta(k+1),-r/2+\delta(k+1)).$$
If $v\notin \cup_{j<k}P_j,$  we define $c_k(v)$ as follows
$$ \left\{\begin{array}{ll}
(-r/2+\delta(k+1), -r/2+\delta(k+1))
 & \mbox{; $vp_k\in E(G), vq_k\in E(G), vs_k\in E(G),$ }\\
(-r+2\delta(k+1)+r(v), -r/2+\delta(k+1)) & \mbox{; $vp_k\notin E(G), vq_k\in E(G), vs_k\in E(G),$ }\\
(-r/2+\delta(k+1), -r+2\delta(k+1)+r(v)) & \mbox{; $vp_k\in E(G), vq_k\notin E(G), vs_k\in E(G),$ }\\
(-r(v),-r(v)) & \mbox{; $vp_k\in E(G), vq_k\in E(G), vs_k\notin E(G),$ }\\
(-r+2\delta(k+1)+r(v),\\
 \hspace{.3in} -r+2\delta(k+1)+r(v))& \mbox{; $vp_k\notin E(G), vq_k\notin E(G), vs_k\in E(G).$ }
\end{array}\right.$$
\end{myenu}

\item \bm{$p_k\notin N(q_k)\nd p_k\notin N^2(q_k).$} Then we define
\begin{myenu}[label=Step \arabic*.]
\item $c_k(p_k):=(-r+2\delta(k+1)-r(p_k),-r+2\delta(k+1)+r(p_k)-\delta).$

\item $c_k(q_k):=(-r+2\delta(k+1)+r(q_k)-\delta ,-r+2\delta(k+1)-r(q_k)).$

\item $c_k(s_k):=(r(s_k),r(s_k)).$

\item Let $np_k\in N(p_k).$
\begin{myenu}[label=Case \arabic*.]
\item If $np_k\in \cup_{j<k}P_j,$ define $$c_k(np_k):=(-r+2\delta(k+1),0).$$

\item If $np_k\in V(G)\sm  \cup_{j<k}P_j,$ then we define
$$c_k(np_k):=\left\{\begin{array}{ll}
(-r+2\delta(k+1),0)  & \mbox{; $np_kq_k\in E(G), np_ks_k\in E(G),$  }\\
(-r(p_k),0) & \mbox{; $np_kq_k\in E(G), np_ks_k\notin E(G),$  }\\
(-r+2\delta(k+1),\\
\hspace{.2in}-r+2\delta(k+1)+r(p_k))   & \mbox{; $np_kq_k\notin E(G), np_ks_k\in E(G).$  }\\
\end{array}\right.$$
\end{myenu}
\item Let $n^2p_k\in N^2(p_k)\setminus \{p_k\}.$
\begin{myenu}[label=Case \arabic*.]
\item If $n^2p_k\in \cup_{j<k}P_j,$ define
$$c_k(n^2p_k):=\left\{\begin{array}{ll}
(-r+2\delta(k+1)+r(p_k),-r/2+\delta (k+1)) & \mbox{; $n^2p_k\in S_u,$ }\\
(-r+2\delta(k+1),0) & \mbox{; $n^2p_k\notin S_u.$ }
\end{array}\right.$$
\item If $n^2p_k\in V(G)\sm  \cup_{j<k}P_j,$ then we define $c_k(n^2p_k) $ as
$$\left\{\begin{array}{ll}
(-r+2\delta(k+1),0) & \mbox{; $n^2p_kq_k\in E(G), n^2p_ks_k\in E(G),$}\\
(-r(p_k),0) & \mbox{; $n^2p_kq_k\in E(G), n^2p_ks_k\notin E(G),$}\\
(-r+2\delta(k+1),-r+2\delta(k+1)+r(p_k)) & \mbox{; $n^2p_kq_k\notin E(G), n^2p_ks_k\in E(G).$}\\
\end{array}\right.$$
\end{myenu}
\item Let $nq_k\in N(q_k).$ (similar to $np_k$)
\begin{myenu}[label=Case \arabic*.]
\item If $nq_k\in \cup_{j<k}P_j,$ define $$c_k(nq_k):=(0,-r+2\delta(k+1)).$$
\item If $nq_k\in V(G)\sm  \cup_{j<k}P_j,$ then we define
$$c_k(nq_k):=\left\{\begin{array}{ll}
(0,-r+2\delta(k+1)) &\mbox{; $nq_kp_k\in E(G), nq_ks_k\in E(G),$}\\
(0,-r(q_k)) &\mbox{; $nq_kp_k\in E(G), nq_ks_k\notin E(G),$}\\
(-r+2\delta(k+1)+r(q_k),\\
\hspace{.3in}-r+2\delta(k+1)) &\mbox{; $nq_kp_k\notin E(G), nq_ks_k\in E(G).$}\\
 \end{array}\right.$$
\end{myenu}
\item Let $n^2q_k\in N^2(q_k)\setminus \{q_k\}.$ (This case is similar to $n^2p_k.$)
\begin{myenu}[label=Case \arabic*.]

\item If $n^2q_k\in \cup_{j<k}P_j,$ then we define
$$c_k(n^2q_k):=\left\{\begin{array}{ll}
(-r/2+\delta (k+1),-r+2\delta(k+1)+r(q_k)) &\mbox{; $n^2q_k\in S_u,$}\\
(0,-r+2\delta(k+1)) &\mbox{; $n^2q_k\notin S_u.$}
 \end{array}\right.$$
 \item If $n^2q_k\in V(G)\sm  \cup_{j<k}P_j,$ then we define
$$c_k(n^2q_k):=\left\{\begin{array}{ll}
(0,-r+2\delta(k+1)) &\mbox{;  $n^2q_kp_k\in E(G), n^2q_ks_k\in E(G),$}\\
(0,-r(q_k)) &\mbox{;  $n^2q_kp_k\in E(G), n^2q_ks_k\notin E(G),$}\\
(-r+2\delta(k+1)+r(q_k),\\
\hspace{.3in}-r+2\delta(k+1)) &\mbox{;  $n^2q_kp_k\notin E(G), n^2q_ks_k\in E(G).$}\\
 \end{array}\right.$$
 \end{myenu}
\item Let $ns_k\in N(s_k).$
\begin{myenu}
\item If $ns_k\in \cup_{j<k}P_j$ then $c_k(ns_k):=(0,0).$
\item If $ns_k\in V(G)\sm  \cup_{j<k}P_j,$ then we define
$$c_k(ns_k):=\left\{\begin{array}{ll}
 (0,0)  &\mbox{;  $ns_kp_k\in E(G), ns_kq_k\in E(G),$}\\
 (-r+2\delta(k+1)+r(s_k),0)    &\mbox{;  $ns_kp_k\notin E(G), ns_kq_k\in E(G),$}\\
    (-r+2\delta(k+1)+r(s_k),\\
\hspace{.3in}-r+2\delta(k+1)+r(s_k)) &\mbox{;  $ns_kp_k\notin E(G), ns_kq_k\notin E(G),$}\\ (0,-r+2\delta(k+1)+r(s_k))  &\mbox{;  $ns_kp_k\in E(G), ns_kq_k\notin E(G).$}
\end{array}\right.$$
\end{myenu}
\item Let $n^2s_k\in N^2(s_k)\setminus \{s_k\}.$
\begin{myenu}[label=Case \arabic*.]
\item If $n^2s_k\in \cup_{j<k}P_j,$ then we define
$$c_k(n^2s_k):=\left\{\begin{array}{ll}
(-r(s_k),-r(s_k))&\mbox{; $n^2s_k\in S_u,$}\\
(0,0)&\mbox{; $n^2s_k\notin S_u.$}
 \end{array}\right.$$

\item If $n^2s_k\in V(G)\sm  \cup_{j<k}P_j,$ then we define
$$c_k(n^2s_k):=
\left\{\begin{array}{ll}
(0,0) &\mbox{; $n^2s_kp_k\in E(G), n^2s_kq_k\in E(G),$}\\
(-r+2\delta(k+1)+r(s_k),0) &\mbox{; $n^2s_kp_k\notin E(G), n^2s_kq_k\in E(G),$}\\
 (0,-r+2\delta(k+1)+r(s_k))&\mbox{; $n^2s_kp_k\in E(G), n^2s_kq_k\notin E(G),$}\\
(-r+2\delta(k+1)+r(s_k),\\
\hspace{.2in} -r+2\delta(k+1)+r(s_k)) &\mbox{; $n^2s_kp_k\notin E(G), n^2s_kq_k\notin E(G).$}
 \end{array}\right.$$
\end{myenu}
\item Let $v\in V(G)\setminus N_k.$ We have the following cases.
If $v\in \cup_{j<k}P_j$ then define
$$c_k(v):=(-r/2+\delta(k+1),-r/2+\delta(k+1)).$$
If $v\notin \cup_{j<k}P_j,$ then we define $c_k(v)$ as follows:
$$\left\{\begin{array}{ll}
(-r/2+\delta(k+1),-r/2+\delta(k+1)) &\mbox{; $vp_k\in E(G), vq_k\in E(G), vs_k\in E(G),$ }\\
  (-r+2\delta(k+1)+r(v),-r/2+\delta(k+1))&\mbox{; $vp_k\notin E(G), vq_k\in E(G), vs_k\in E(G),$ }\\
  (-r/2+\delta(k+1),-r+2\delta(k+1)+r(v)) &\mbox{; $vp_k\in E(G), vq_k\notin E(G), vs_k\in E(G),$ }\\
 (-r(v),-r(v))   &\mbox{; $vp_k\in E(G), vq_k\in E(G), vs_k\notin E(G),$ }\\
 (-r+2\delta(k+1)+r(v),\\
\hspace{.2in}-r+2\delta(k+1)+r(v))    &\mbox{; $vp_k\notin E(G), vq_k\notin E(G), vs_k\in E(G).$ }
\end{array}\right.$$

 \end{myenu}
\end{myenu}

\subsection{Proof of (\ref{LemEq1})}

Let $P_k:=\{p_k,q_k,s_k\}$ be  a triplet with   exactly one edge among its vertices, picked up through the recursion of    step   \ref{CharVerS24} of Algorithm \ref{alg41}.
Without loss of generality, assume that $p_kq_k\in E(G).$

Let $u\in V(G) \nd nu\in N(u).$ We need to establish that
\begin{equation*}
 |c_k(u)-c_k(nu)|  \leq r(u).
\end{equation*}

%
Let $u \in P_i\nd nu \in P_j.$
Then by definition $m(u)\leq \min\{i,j\}$
and thence $$r(u)=r-2\del m(u) \geq r-2\del \min\{i,j\}.$$

We now study  case by case as follows.

\begin{myenu}[label=Case \arabic*.]
\item \bm{$p_k\in N(q_k)\mbox{ or }p_k\in N^2(q_k).$}
\begin{myenu}
\item If $u\in N(p_k),$ then
\begin{align*}
c_k(u)\in     \big{\{}(-r+&2\delta(k+1)-\delta,-r+2\delta(k+1)-r(u)),\\
\hspace{.2in} &(-r+2\delta(k+1),-r+2\delta(k+1)), (-r(u),-r(u)) \big{\}} \\
\nd  c_k(nu) \in \big{\{}(-r+&2\delta(k+1)-r(u),-r+2\delta(k+1)),\\
 \hspace{.2in}  &(-r+2\delta(k+1),-r+2\delta(k+1)), (-r(u),-r(u)) \big{\}}.
 \end{align*}

\item Similarly, if $u\in N(q_k),$ then above equation holds true.

\item If $u\in N(s_k),$ then
\begin{align*}
c_k(u)\in \big{\{}(0&,0),(-r+2\delta(k+1)+r(u),0),(0,-r+2\delta(k+1)+r(u)),\\
\hst &(-r+2\delta(k+1)+r(u),-r+2\delta(k+1)+r(u))\big{\}}.\\
\nd c_k(nu)\in  \big{\{}(r&(u),r(u)),(-r+2\delta(k+1)+r(u)-\delta ,-r(u)),(0,0),\\
\hst &(-r+2\delta(k+1)+r(u),0), (0,-r+2\delta(k+1)+r(u)),\\
\hst &(-r+2\delta(k+1)+r(u),-r+2\delta(k+1)+r(u))\big{\}}.
\end{align*}
The required equation does not hold, only for the case when $c_k(u)=(0,-r+2\delta(k+1)+r(u))$ and $c_k(nu)=(-r+2\delta(k+1)+r(u)-\delta ,-r(u)),$ (consider the  second co-ordinates  of both of these).

However, such a situation is impossible, as
when $nu$ has this co ordinate, then it is an exterior vertex of some star and  the center of that star is already picked earlier due to our algorithm. Hence such co-ordinates for $u$ and $nu$ are not possible together for this module.
Moreover, whenever $c_k(nu)$ has these co-ordinates, then $c_k(u)=(0,0).$
\item\label{1.1.4.} If $u\notin \{p_k,q_k,t_k\}\cup N(s_k)\cup N^2(s_k)=Z,$ (say) then $nu\notin Z.$
\\Therefore, each of $c_k(u)[0],c_k(u)[1],c_k(nu)[0]\nd c_k(nu)[1]$ belongs to set
$$\big{\{}-r/2+\delta(k+1),-r+2\delta(k+1)+r(u),-r(u)\big{\}}.$$
Again the required equation does not hold, only for the case when $c_k(u)-c_k(nu)=-r+2\delta(k+1)+r(u)-(-r(u)).$
Note that we are picking both $u$ and $nu$ after $P_k.$
Consider the following cases.
\begin{myenu}
  \item If $u$ and $nu$ are vertices of a matching edge, then
$m(u)$ is greater than $k.$ Hence $r(u)\leq r-2\delta (k+1).$ Hence the inequality holds.
\item If $u$ and $nu$ are vertices of triangle, then by our definitions, again \wv  $r(u)\leq r-2\delta (k+1).$ Hence no issues.
\item If $u$ and $ nu$ are  vertices of a star, then due to our characterizing vertices algorithm one of $u$ or $nu$ comes before $P_k.$ Hence this case fails.
\end{myenu}
\end{myenu}
\item \bm{$p_k\notin N(q_k)\nd p_k\notin N^2(q_k).$}
\begin{myenu}
\item If $u\in N(p_k),$ then
\begin{align*}
c_k(u)\in \bigg{\{}(-r&+2\delta(k+1),0),(-r(u),0),\\
\hst &(-r+2\delta(k+1),-r+2\delta(k+1)+r(u))\bigg{\}}.\\
c_k(nu)\in   \bigg{\{}(-r&+2\delta(k+1)-r(u),-r+2\delta(k+1)+r(u)-\delta),\\
\hst &(-r+2\delta(k+1)+r(u),-r/2+\delta (k+1)),(-r(u),0),\\
\hst &(-r+2\delta(k+1),0),(-r+2\delta(k+1),-r+2\delta(k+1)+r(u))\bigg{\}}
\end{align*}
Note that, when $c_k(nu)=(-r+2\delta(k+1)+r(u),-r/2+\delta (k+1)),$ then $$c_k(u)=(-r+2\delta(k+1),0).$$

\item Similarly, if $u\in N(q_k),$ then above equation holds true.

\item If $u\in N(s_k),$ then
\begin{align*}
c_k(u)\in \bigg{\{}(0,&0),(-r+2\delta(k+1)+r(u),0),(0,-r+2\delta(k+1)+r(u)),\\
\hst &(-r+2\delta(k+1)+r(u),-r+2\delta(k+1)+r(u))\bigg{\}}.\\
c_k(nu)\in   \bigg{\{}(r&(u),r(u)),(-r(u),-r(u)),(0,0),(-r+2\delta(k+1)+r(u),0),(0,-r\\
\hst &+2\delta(k+1)+r(u)),(-r+2\delta(k+1)+r(u),-r+2\delta(k+1)+r(u))\bigg{\}}.
\end{align*}
Note that, when $c_k(nu)=(-r(u),-r(u)),$ then $c_k(u)=(0,0).$

\item  If $u\notin \{p_k,q_k,t_k\}\cup N(s_k)\cup N^2(s_k),$  then this case is similar to \ref{1.1.4.}  above.
\end{myenu}
\end{myenu}

\subsection{Proof of (\ref{LemEq2})}

Let  $u\in P_k$  and  $v\in \displaystyle\bigcup_{l\leq k} P_k\sm \{u\}.$
 We need to  prove that
\begin{equation*}
 |c_k(v)-c_k(u)|  \geq \max\{r(u),r(v)\}.
\end{equation*}

%
We now study this case by case as follows.
%
\begin{myenu}[label=Case \arabic*.]
\item \bm{$p_k\in N(q_k)\mbox{ or }p_k\in N^2(q_k).$} Then $c_k(u)$ belongs to the set
\begin{align*}
 {\{}(-r&+2\delta (k+1)-r(u),-r+2\delta (k+1)),(r(u),r(u)),\\
& (-r+2\delta (k+1)-\delta ,-r+2\delta (k+1)-r(u)) {\}}.
\end{align*}
\begin{myenu}
\item If $c_k(u)= (-r+2\delta (k+1)-r(u),-r+2\delta (k+1)),$ then $c_k(v)$ belongs to the set
\begin{align*}
 \bigg{\{}(-&r+2\delta (k+1)-\delta ,-r+2\delta (k+1)-r(v)),(-r+2\delta (k+1),-r+2\delta (k+1)),\\
& (r(v),r(v)),(0,0),(-r+2\delta (k+1)+r(v)-\delta ,-r(v)), \big(\frac{r}{2}+\delta (k+1),\frac{r}{2}+\delta (k+1) \big) \bigg{\}}.
\end{align*}
\item If $c_k(u)= (-r+2\delta (k+1)-\delta ,-r+2\delta (k+1)-r(u)),$ then $c_k(v)$ belongs to
\begin{align*}
\bigg{\{}(-r&+2\delta (k+1)-r(v),-r+2\delta (k+1)),
(-r+2\delta (k+1),-r+2\delta (k+1)),\\
&(r(v),r(v)),(0,0),(-r+2\delta (k+1)+r(v)-\delta ,-r(v)),\\
&\bigg(-\frac{r}{2}+\delta (k+1),-\frac{r}{2}+\delta (k+1)\bigg) \bigg{\}}.
\end{align*}
\item  If $c_k(u)= (r(u),r(u)),$ then $c_k(v)$ belongs to
\begin{align*}
 \bigg{\{}(-r&+2\delta (k+1)-r(v),-r+2\delta (k+1)),(-r+2\delta (k+1)-\delta ,-r+2\delta (k+1)-r(v)),\\
&(-r+2\delta (k+1),-r+2\delta (k+1)), (-r+2\delta (k+1)+r(v)-\delta ,-r(v)),\\
&(0,0),\bigg(-\frac{r}{2}+\delta (k+1),-\frac{r}{2}+\delta (k+1)\bigg) \bigg{\}}.
\end{align*}
\end{myenu}
\item \bm{$p_k\not\in N(q_k)\nd p_k\not\in N^2(q_1).$} Then
\begin{align*}
c_k(u)\in   {\{}(-r&+2\delta (k+1)-r(u),-r+2\delta (k+1)+r(u)-\delta ),\\
&(-r+2\delta (k+1)+r(u)-\delta ,-r+2\delta (k+1)-r(u)),(r(u),r(u)) {\}}.
 \end{align*}
\begin{myenu}
\item If $c_k(u)= (-r+2\delta (k+1)-r(u),-r+2\delta (k+1)+r(u)-\delta ),$ then
$c_k(v)$ belongs to
\begin{align*}
  \bigg{\{}(-r&+2\delta (k+1)+r(v)-\delta ,-r+2\delta (k+1)-r(v)),(r(v),r(v)),(-r+2\delta (k+1),0),\\
&\bigg(-r+2\delta (k+1)+r(v),\frac{r}{2}+\delta (k+1)\bigg),(-r+2\delta (k+1),0),(0,-r+2\delta (k+1)),\\
& \bigg(\frac{r}{2}+\delta (k+1),-r+2\delta (k+1)+r(v)\bigg),(0,-r+2\delta (k+1)),(0,0),(-r(v),-r(v)),\\
& \bigg(-\frac{r}{2}+\delta (k+1),-\frac{r}{2}+\delta (k+1)\bigg) \bigg{\}}.
\end{align*}
Recall that in case
 $c_k(v)=(-r(v),-r(v))$ and  $$c_k(u)= (-r+2\delta (k+1)-r(u),-r+2\delta (k+1)+r(u)-\delta ),$$
\label{p44ineq5}the calculations for the required inequality are already done on page \pageref{pageconj14}. The other  calculations are not that tricky.

%
%
 \item If $c_k(u)=(-r+2\delta (k+1)+r(u)-\delta ,-r+2\delta (k+1)-r(u)),$ then $c_k(v)$ belongs to
\begin{align*}
 \bigg{\{}(-r&+2\delta (k+1)-r(v),-r+2\delta (k+1)+r(v)-\delta ),(r(v),r(v)),(-r+2\delta (k+1),0),\\
& \bigg(-r+2\delta (k+1)+r(v),-\frac{r}{2}+\delta (k+1)\bigg),(-r+2\delta (k+1),0),(0,-r+2\delta (k+1)),\\
& \bigg(-\frac{r}{2}+\delta (k+1),-r+2\delta (k+1)+r(v)\bigg),(0,-r+2\delta (k+1)),(0,0),(-r(v),-r(v)),\\
& \bigg(-\frac{r}{2}+\delta (k+1),-\frac{r}{2}+\delta (k+1)\bigg) \bigg{\}}.
\end{align*}
\item  If $c_k(u)=(r(u),r(u)),$ then $c_k(v)$ belongs to the set
\begin{align*}
 \bigg{\{}(&-r+2\delta (k+1)-r(v),-r+2\delta (k+1)+r(v)-\delta ),(-r+2\delta (k+1),0),(0,0),\\
 &(-r+2\delta (k+1)+r(v)-\delta, -r+2\delta (k+1)-r(v)),(-r+2\delta (k+1),0),(-r(v),-r(v)),\\
 &\bigg(-r+2\delta (k+1)+r(v),-\frac{r}{2}+\delta (k+1)\bigg), \bigg(-\frac{r}{2}+\delta (k+1),-r+2\delta (k+1)+r(v)\bigg),\\
&(0,-r+2\delta (k+1)), (0,-r+2\delta (k+1)),\bigg(-\frac{r}{2}+\delta (k+1),-\frac{r}{2}+\delta (k+1)\bigg)\bigg{\}}.
\end{align*}
\end{myenu}
\end{myenu}


%

%

\subsection{Proof of (\ref{LemEq3})}

Let  $u\in P_k\nd  v\in \displaystyle\bigcup_{l\leq k} P_l \sm \{u\}$  such that $uv\not\in E(G)$   and $\{u,v\}\subseteq S_w \fos w\in V(G).$
\\ We need to  prove that
$$|c_k(v)-c_k(u)|  \geq r(v)+r(u).$$

%
We now study this case by case as follows.

Since $\{u,v\}\subset S_w, \mbox{ where } w\notin \{u,v\}.,$
\wo  $v\in N^2(u)\setminus \{u\}.$

\begin{myenu}[label=Case \arabic*.]
\item \bm{$p_k\in N(q_k)\mbox{ or }p_k\in N^2(q_k)\setminus \{q_k\}.$}
\begin{myenu}
\item \bm{$u=p_k.$} Then $p_k\in S_w.$ Since $p_kq_k\in E(G),$ \wo $q_k\in N(p_k)\subset N^2(w)=\{w\}.$ Hence $q_k=w,$ a central vertex. This is impossible, as all central vertices of our stars had been chosen before step        \ref{CharVerS24} of Algorithm \ref{alg41}, when $q_k$ was chosen. Hence $u\neq p_k.$
 \item \bm{$u=q_k.$}  As above, this case is also impossible.
\item \bm{$u=s_k.$} Then
$$c_k(u)=(r(u),r(u))\nd c_k(v)=(-r+2\delta (k+1)+r(v),-r(v)).$$
\end{myenu}
\item \bm{$p_k\not\in N(q_k)\mbox{ and }p_k\not\in N^2(q_k)\setminus \{q_k\}.$}
\begin{myenu}
\item \bm{$u=p_k.$} Then
\begin{align*}
   c_k(u)=&(-r+2\delta (k+1)-r(u),-r+2\delta (k+1)+r(u)-\delta )\\
\nd c_k(v)=&(-r+2\delta (k+1)+r(v),-r/2+\delta (k+1)).
\end{align*}
\item \bm{$u=q_k.$} This case is  similar to above case, as \wv
\begin{align*}
   c_k(u)=&(-r+2\delta (k+1)+r(u)-\delta ,-r+2\delta (k+1)-r(u))\\
\nd c_k(v)=&(-r/2+\delta (k+1),-r+2\delta (k+1)+r(v)).
\end{align*}
\item \bm{$u=s_k.$} Then
$$c_k(u)=(r(u),r(u))\nd c_k(v)=(-r(v),-r(v)).$$
\end{myenu}
\end{myenu}

 \subsection{Proof of (\ref{LemEq4})}

Let  $u\in P_k$  and  $v\in \displaystyle\bigcup_{l\geq k} P_l\sm \{u\} $ be such that $uv\not\in E(G)$   and $\{u,v\}\not \subseteq S_w \foa w\in V(G).$
\\We need to  prove that
$$|c_k(v)-c_k(u)|  \geq r(v)+r(u).$$
We now study this case by case as follows.
\begin{myenu}[label=Case \arabic*.]
\item \bm{$p_k\in N(q_k)\mbox{ or }p_k\in N^2(q_k)\setminus \{q_k\}.$}
\begin{myenu}
\item \bm{$u=p_k.$} Then $$c_k(u)=(-r+2\delta(k+1)-r(u),-r+2\delta(k+1))$$
 and $c_k(v)$ belongs to the set
\begin{align*}
 \big{\{}(r(v),&r(v)),(-r+2\delta(k+1)+r(v),0), (-r+2\delta(k+1)+r(v),-r/2+\delta(k+1)),\\
&  (-r+2\delta(k+1)+r(v),-r+2\delta(k+1)+r(v)) \big{\}}.
\end{align*}
\item \bm{$u=q_k.$} This case is similar to above case.
\item \bm{$u=s_k.$} Then $c_k(u)=(r(u),r(u))$ and $c_k(v)$ belongs to
\begin{align*}
 \big{\{}(-r&+2\delta(k+1)-r(v),-r+2\delta(k+1)), (-r(v),-r(v))\\
&(-r+2\delta(k+1)-\delta ,-r+2\delta(k+1)-r(v)), \big{\}}.
\end{align*}
\end{myenu}
\item \bm{$p_k\not\in N(q_k)\nd p_k\not \in N^2(q_k)\setminus \{q_k\}.$}
\begin{myenu}
\item \bm{$u=p_k.$} Then
$$c_k(u)=(-r+2\delta(k+1)-r(u),-r+2\delta(k+1)+r(u)-\delta)$$
and $c_k(v)$ belongs to the set
\begin{align*}
  \big{\{}(r&(v),r(v)),  (-r+2\delta (k+1)+r(v),0),
  (-r+2\delta(k+1)+r(v),  -r/2+\delta(k+1)),\\
&  (-r+2\delta(k+1)+r(v),-r+2\delta(k+1)),
  (-r+2\delta (k+1)+r(v),-r+2\delta (k+1)+r(v))
  \big{\}}.
\end{align*}
\item \bm{$u=q_k.$} This case is similar to the previous case. \Nt
\\$c_k(u)=(-r+2\delta(k+1)+r(u)-\delta,-r+2\delta(k+1)-r(u))\nd c_k(v)$ belongs to the set
\begin{align*}
   \big{\{}(r&(v),r(v)),
   (0,-r+2\delta(k+1)+r(v)),
   (-r+2\delta(k+1),-r+2\delta(k+1)+r(v)),\\
&   (-r/2+\delta(k+1),-r+2\delta(k+1)+r(v)),
    (-r+2\delta(k+1)+r(v), -r+2\delta(k+1)+r(v))
  \big{\}}.
\end{align*}
\item \bm{$u=s_k.$} Then $c_k(u)=(r(u),r(u))\nd c_k(v)$ belongs to the set
\begin{align*}
  \big{\{}
  (-r&+2\delta(k+1)-r(v),-r+2\delta(k+1)+r(v)-\delta),(-r(v),0),  (0,-r(v)),(-r(v),-r(v)),\\
&  (-r+2\delta(k+1)+r(v)-\delta, -r+2\delta(k+1)-r(v))\big{\}}.
\end{align*}
\end{myenu}
\end{myenu}

\subsection{Proof of (\ref{LemEq5})}

Let   $u\in V(G)\nd v\in V(G) \sm \{u\}$ such that $uv\in E(G).$  We will  prove that
$$|c_k(u)-c_k(v)|  < r(u)+r(v).$$

Then there exist $i$ and $j$ \st
$u \in P_i\nd v \in P_j.$
\Csq $$r(u) \geq r-2\delta i \mbox{ and } r(v) \geq r-2\delta j.$$
We now study this case by case as follows.

\begin{myenu}[label=Case \arabic*.]
\item \bm{$p_k\in N(q_k)\mbox{ or }p_k\in N^2(q_k)\setminus \{q_k\}.$}
\begin{myenu}
\item \bm{$k\notin \{i,j\}.$} Let $Q:=(-r/2+\delta (k+1),-r/2+\delta (k+1))\nd A(a)$ denotes the set
\begin{align*}
   \big{\{}
(-r&+2\delta(k+1),-r+2\delta(k+1)),(-r(a),-r(a)),(0,0),(-r+2\delta (k+1)+r(a),0),\\
&(0,-r+2\delta (k+1)+r(a)),(-r+2\delta (k+1)+r(a),-r+2\delta (k+1)+r(a)),\\
&(-r+2\delta (k+1)+r(a)-\delta ,-r(a)),(-r+2\delta (k+1)+r(a),-r/2+\delta (k+1)), \\
& (-r/2+\delta (k+1),-r/2+\delta (k+1)),(-r/2+\delta (k+1),-r+2\delta (k+1)+r(a))\big\}.
\end{align*}
\Nt we have $c_k(u)\in A(u)\nd c_k(v)\in A(v).$
Observe that
$$|c_k(u)-Q|<r(u)\nd |c_k(v)-Q|<r(v).$$
 Applying the triangle inequality, \wo the required inequality.

\item \bm{$k\in \{i,j\}.$} Without loss of generality, we assume that $k=i.$ Then   $c_k(u)$ belongs to
\begin{align*}
    \{(&-r+2\delta (k+1)-r(u),-r+2\delta (k+1)),(r(u),r(u))\\
    &(-r+2\delta (k+1)-\delta,    -r+2\delta (k+1)-r(u))\}.
    \end{align*}
If $c_k(u)=(-r+2\delta (k+1)-r(u),-r+2\delta (k+1)),$ then $c_k(v)$ belongs to the set
 \begin{align*}
   \big{\{}
   (&-r+2\delta (k+1)-\delta ,-r+2\delta (k+1)-r(v)), (-r+2\delta(k+1),-r+2\delta(k+1)),\\
&  (-r(v),-r(v)),(0,0),(0,-r+2\delta (k+1)+r(v)),(-r+2\delta (k+1)+r(v)-\delta ,-r(v)), \\
&   (-r/2+\delta (k+1),-r/2+\delta (k+1)),   (-r/2+\delta (k+1),-r+2\delta (k+1)+r(v))
   \big{\}}.
\end{align*}

Finally, if $$c_k(u)=(-r+2\delta (k+1)-\delta ,-r+2\delta (k+1)-r(u)),$$
 then $c_k(v)$ belongs to the set
\begin{align*}
   \big{\{}
  & (-r+2\delta (k+1)-r(v),-r+2\delta (k+1)), (-r+2\delta(k+1),-r+2\delta(k+1)),\\
& (-r(v),-r(v)), (-r+2\delta (k+1)+r(v),0),(-r+2\delta (k+1)+r(v)-\delta ,-r(v)),   \\
 &(-r+2\delta (k+1)+r(v),0),(0,0),(-r/2+\delta (k+1),-r/2+\delta (k+1)),\\
 &(-r+2\delta (k+1)+r(v),-r/2+\delta (k+1))\big{\}}.
\end{align*}
If $c_k(u)=(r(u),r(u)),$ then $c_k(v)$ belongs to the set
 \begin{align*}
   \big{\{}
 (&-r+2\delta(k+1),-r+2\delta(k+1)),  (-r+2\delta (k+1)+r(v),0),(0,-r+2\delta (k+1)+r(v)),\\
&(-r+2\delta (k+1)+r(v),-r+2\delta (k+1)+r(v)),(-r/2+\delta (k+1),-r/2+\delta (k+1)),(0,0), \\
&(-r+2\delta (k+1)+r(v),-r/2+\delta (k+1)),(-r/2+\delta (k+1),-r+2\delta (k+1)+r(v)) \big{\}}.
\end{align*}
Let $Q:=(\delta ,\delta ).$ Observe that
$$|c_k(u)-Q|<r(u)\nd |c_k(v)-Q|<r(v).$$
Hence $|c_k(u)-c_k(v)|<r(u)+r(v).$
\end{myenu}
\item \bm{$p_k\not\in N(q_k)\mbox{ and }p_k\not\in N^2(q_k)\setminus \{q_k\}.$}
\begin{myenu}
\item \bm{$k\notin \{i,j\}.$} Let $Q:=(-r/2+\delta (k+1),-r/2+\delta (k+1))\nd A(a)$ denotes the set
\begin{align*}
  \big{\{}
 (-r&+2\delta(k+1),0), (-r(a),0),(-r+2\delta(k+1),-r+2\delta(k+1)+r(a)),(0,-r(a)),\\
& (-r+2\delta(k+1)+r(a),-r/2+\delta (k+1)),  (-r+2\delta(k+1)+r(a),-r+2\delta(k+1)),\\
&(-r/2+\delta (k+1),-r+2\delta(k+1)+r(a)),   (-r(a),-r(a)), (-r+2\delta (k+1)+r(a),0),\\
&(0,0), (-r+2\delta (k+1)+r(a),-r+2\delta (k+1)+r(a)),(0,-r+2\delta (k+1)+r(a)),\\
&(-r/2+\delta (k+1),-r/2+\delta (k+1)),(-r+2\delta (k+1)+r(a),-r/2+\delta (k+1)), \\
&(-r/2+\delta (k+1),-r+2\delta (k+1)+r(a)), (0,-r+2\delta(k+1)) \big{\}}.
\end{align*}
Then $c_k(u)\in A(u)\nd c_k(v)\in A(v).$ Also observe that
$$|c_k(u)-Q|<r(u)\nd |c_k(v)-Q|<r(v).$$
 Hence the required inequality is satisfied.
\item \bm{$k\in \{i,j\}.$} Without loss of generality, we assume that $k=i.$
\begin{align*}
c_k(u)\in & \big{\{}(-r+2\delta (k+1)-r(u),-r+2\delta (k+1)+r(u)-\delta ),\\
&(-r+2\delta (k+1)+r(u)-\delta ,-r+2\delta (k+1)-r(u)),(r(u),r(u))\big{\}}.
\end{align*}
If $c_k(u)=(-r+2\delta (k+1)-r(u),-r+2\delta (k+1)+r(u)-\delta ),$ then $c_k(v)$ belongs to
 \begin{align*}
  \big{\{}
  (-r&+2\delta (k+1)+r(v)-\delta ,-r+2\delta (k+1)-r(v)),(-r+2\delta(k+1),0), (-r(v),0),\\
&  (-r+2\delta(k+1),-r+2\delta(k+1)+r(v)),(-r+2\delta(k+1),-r+2\delta(k+1)+r(v)),\\
&  (0,-r+2\delta(k+1)), (0,-r(v)), (-r/2+\delta(k+1),-r+2\delta(k+1)+r(v)),   (0,0),\\
 & (0,-r+2\delta (k+1)+r(v)),(-r(v),-r(v)), (-r/2+\delta (k+1),-r/2+\delta (k+1)), \\
 & (-r/2+\delta (k+1),-r+2\delta (k+1)+r(v)) \big{\}}.
\end{align*}
If $$c_k(u)=(-r+2\delta (k+1)+r(u)-\delta ,-r+2\delta (k+1)-r(u)),$$ then $c_k(v)$ belongs to
 \begin{align*}
   \big{\{}
 (-r&+2\delta (k+1)-r(v),-r+2\delta (k+1)+r(v)-\delta),(-r+2\delta(k+1),0),\\
 &  (-r(v),0),(-r+2\delta(k+1)+r(v),-r/2+\delta k+1)), (0,-r+2\delta(k+1)),\\
&(0,-r(v)),(-r+2\delta(k+1)+r(v),-r+2\delta(k+1)), (0,0),(-r(v),-r(v)),\\
& (-r+2\delta (k+1)+r(v),0), (-r/2+\delta (k+1),-r/2+\delta (k+1)),\\
& (-r+2\delta (k+1)+r(v),-r/2+\delta (k+1))
\big{\}}.
\end{align*}
%
Similarly, if $c_k(u)=(r(u),r(u)),$ then $c_k(v)$ belongs to
 \begin{align*}
\big{\{}
(-r&+2\delta(k+1),0), (-r+2\delta(k+1),-r+2\delta(k+1)+r(v)), (0,-r+2\delta(k+1)), \\ &(-r+2\delta(k+1)+r(v),-r+2\delta(k+1)),(-r/2+\delta (k+1),-r+2\delta(k+1)+r(v)),\\
& (0,0),(-r+2\delta (k+1)+r(v),0),(-r+2\delta (k+1)+r(v),-r+2\delta (k+1)+r(v)),\\
& (0,-r+2\delta (k+1)+r(v)),  (-r/2+\delta (k+1),-r/2+\delta (k+1)), (-r+2\delta (k+1)+r(v),\\
&-r/2+\delta (k+1)),(-r/2+\delta (k+1),-r+2\delta (k+1)+r(v)),(-r+2\delta (k+1)+r(v),\\
&-r+2\delta (k+1)+r(v)),(-r+2\delta(k+1)+r(v),-r/2+\delta (k+1))\big{\}}.
\end{align*}
Let $Q:=(\delta ,\delta ).$ \Nt $|c_k(u)-Q|<r(u)\nd |c_k(v)-Q|<r(v).$
Hence
$$|c_k(u)-c_k(v)|\leq |c_k(u)-Q|+|Q-c_k(u)|<r(u)+r(v).$$
As usual, the triangle inequality implies the required inequality.
\end{myenu}

\end{myenu}

\section{\bf Module VIII: The Super Triplets}\label{sectionM8}


%
%

Finally, we consider the case of  triplets $P_i,$
picked up either till step (\ref{CharVerS17})  or  during the recursion of step (\ref{CharVerS22}) of Algorithm \ref{alg41}.
Recall that  $P_i$ consists of independent vertices and  no two vertices out of these  belong to a single star.
These triplets will be called the \textit{super triplets}, throughout this module.

The following are these cases, of super triplets, which we further examine as per our own requirement.

\begin{myenu}[label=Case \Roman*.] \item Consider the  $i^{th}$ triplet $P_i,$ picked up before step (\ref{CharVerS17}) of Algorithm \ref{alg41}.
\begin{myenu}[label=Subcase \roman*.]
\item There exists $ v'\in P_i$ such that $v'\in S_u,$ for some $S_u.$ Let $b:=v'\nd a, c$ be the other two vertices in $P_i.$ This case is a subset of \textit{super triplet's case} presented in the subsequent discussion. If there are more than one such $v',$ then we call any one of them as $b$ and call other two as $a\nd c.$

\item If there does not exist any $v'\in P_i$ such that $v'\in S_u,$ then denote any permutation of $P_i$ as $\{a,b,c\}.$

\end{myenu}

\item Consider $\{p_i,q_i,s_i\},$ when there is \textit{no edge among them}.
Then denote any permutation of $P_i$ as $\{a,b,c\},$
say $a:=p_i,b:=q_i\nd c:=s_i.$

\end{myenu}

\subsection{Assigning Dimensions for Super Triplets}\label{assdimsec}


We now demonstrate as how to assign dimensions to   the vertices picked up in Algorithm \ref{alg41} and how this works in embedding our given graph into some Euclidean space as a sphere of influence graph. We present the particular example of the super triplet's case, which covers most of the cases.

Consider the case of assigning two dimensions (ordinates) corresponding to the super triplet  $P_i=\{a, b, c\},$ where $a,b,c$ do not have an edge among them and no two of $\{a,b,c\}$ are in same $S_u.$ We define $c_i(v)$ for all $v\in V(G),$ as follows.


\begin{myenu}[label= Step \arabic*.]
\item Define
$c_i(a):=(-r+2\delta(i+1)-r(a), r(a)).$
\item Define
$c_i(b):=(\delta+r(b), \delta+r(b)).$
\item Define
$c_i(c):=(r(c),-r+2\delta(i+1)-r(c)).$
\item
Let $na\in N(a)$ be arbitrary. Define $c_i(na),$ as follows.
\begin{myenu} [label= Case  \arabic*.]
\item If $na\in \cup_{j<i}P_j,$ then we define $c_i(na):=(-r+2\delta(i+1),0).$
\item If $na\notin \cup_{j<i}P_j,$   define $c_i(na)$ as follows.
$$\left\{\begin{array}{ll}
(-r(a),0) & \mbox{; $nab\notin E(G)\nd nac \notin E(G),$ }\\
(-r+2\delta (i+1),0) &  \mbox{; $nab\in E(G)\nd nac\in E(G),$ }\\
(-r(a)+\delta ,0) & \mbox{; $nab\notin E(G)\nd nac\in E(G),$ }\\
(-r+2\delta(i+1),-r+2\delta(i+1)+r(a))   & \mbox{; $nab\in E(G)\nd nac\notin E(G).$} \end{array}\right.$$
\end{myenu}

\item
Let $nb\in N(b)$ be arbitrary. Define $c_i(nb),$ as follows.
\begin{myenu} [label= Case  \arabic*.]
\item If $nb\in \cup_{j<i}P_j,$   define $c_i(nb):=(\delta ,\delta ).$
\item If $nb\notin \cup_{j<i}P_j,$   define $c_i(nb)$ as follows.
$$\left\{\begin{array}{ll}
(\delta, \delta) & \mbox{; $nba\in E(G)\nd nbc\in E(G),$ }\\
(\delta, -r+2\delta(i+1)+r(b)) &  \mbox{; $nba\in E(G)\nd nbc\notin E(G),$ }\\
(-r+2\delta(i+1)+r(b), \delta) & \mbox{; $nba\notin E(G)\nd nbc\in E(G),$} \\
(-r+2\delta(i+1)+r(b), -r+2\delta(i+1)+r(b)) & \mbox{; $nba\notin E(G)\nd nbc\notin E(G)$}\\ & \mbox{ and $nb$ is not a central vertex.}
\end{array}\right.$$
\end{myenu}

\item Define   {$c_i$ on $nc\in N(c),$} similar to that of $na\in N(a),$ as follows.
\begin{myenu} [label= Case  \arabic*.]
\item If $nc\in \cup_{j<i}P_j,$ define $c_i(nc)=(0,-r+2\delta(i+1)).$
\item If $nc\notin \cup_{j<i}P_j,$   we define $c_i(nc)$ as follows.
$$ \left\{\begin{array}{ll}
(0,-r(c)) & \mbox{; $nca\notin E(G)\nd ncb\notin E(G),$ }\\
(0,-r+2\delta (i+1)) &  \mbox{; $nca\in E(G)\nd ncb\in E(G),$ }\\
(0,-r(c)+\delta ) & \mbox{; $nca\in E(G)\nd ncb\notin E(G),$ }\\
(-r+2\delta(i+1)+r(c),-r+2\delta(i+1))   & \mbox{; $nca\notin E(G)\nd ncb\in E(G).$} \end{array}\right.$$
\end{myenu}

\item
Let $n^2a\in N^2_0(a)$  be arbitrary,
$na\in N(a)\nd c_i(na)=(na_x, na_y).$
Set $c_i(n^2a),$ as follows.

\begin{myenu} [label= Case  \arabic*.]
\item If $n^2a\in \cup_{j<i}P_j$ and for some  $u,n^2a\in S_u.$
 We define
$$c_i(n^2a):=\left\{\begin{array}{ll}
(-r+2\delta (i+1), -r(a))     & \mbox{if $na_x\not =-r+2\delta (i+1),$  }\\
(-r+2\delta(i+1)+r(a),-r/2+\delta(i+1)) & \mbox{if $na_x =-r+2\delta (i+1).$} \end{array}\right.$$

\item If $n^2a\in \cup_{j<i}P_j$ and for all  $u,n^2a\notin S_u.$
 We define $c_i(n^2a):=(-r+2\delta (i+1),0).$
\item If $n^2a\notin \cup_{j<i}P_j,$ we define $c_i(n^2a)$ as follows.
$$ \left\{\begin{array}{ll}
(-r(a),0)      & \mbox{; $n^2ab\notin E(G)\nd n^2ac\notin E(G),$ }\\
(-r+2\delta (i+1),0) & \mbox{; $n^2ab\in E(G)\nd n^2ac\in E(G),$ }\\
(-r(a)+\delta ,0) & \mbox{; $n^2ab\notin E(G)\nd n^2ac\in E(G),$ }\\
(-r+2\delta(i+1),-r+2\delta(i+1)+r(a))   & \mbox{; $n^2ab\in E(G)\nd n^2ac\notin E(G).$ } \end{array}\right.$$
 \end{myenu}

\item
Let $n^2b\in N^2_0(b)$  be arbitrary,
 $nb\in N(b)$ and  $c_i(nb)=(nb_x,nb_y).$
Set $c_i(n^2b),$ as follows.

\begin{myenu} [label= Case  \arabic*.]
\item If $n^2b\in \cup_{j<i}P_j$ and for some $u, n^2b\in S_u,$
 we define $c_i(n^2b)$ as follows.
$$ \left\{\begin{array}{ll}
(\delta -r(b), \delta -r(b)) & \mbox{; $c_i(nb)=(\delta , \delta ),$ }\\
(\delta -r(b), -r+2\delta(i+1)) &  \mbox{; $c_i(nb)=(\delta , -r+2\delta(i+1)+r(b)),$ }\\
(-r+2\delta(i+1), \delta -r(b)) & \mbox{; $c_i(nb)=(-r+2\delta(i+1)+r(b), \delta ).$} \end{array}\right.$$

\item If $n^2b\in \cup_{j<i}P_j$ and for all $u, n^2b\notin S_u,$ we define $c_i(n^2b):=(\delta, \delta).$
\item If $n^2b\notin \cup_{j<i}P_j,$ we define $c_i(n^2b)$ as follows.
$$\left\{\begin{array}{ll}
(\delta, \delta) & \mbox{; $n^2ba\in E(G)\nd n^2bc\in E(G),$ }\\
(\delta, -r+2\delta(i+1)+r(b)) &  \mbox{; $n^2ba\in E(G)\nd n^2bc\notin E(G),$ }\\
(-r+2\delta(i+1)+r(b), \delta) & \mbox{; $n^2ba\notin E(G)\nd n^2bc\in E(G),$} \\
(-r+2\delta(i+1)+r(b),-r+2\delta(i+1)+r(b)) & \mbox{; $n^2ba\notin E(G)\nd n^2bc\notin E(G).$}\end{array}\right.$$
\end{myenu}

\item
Let $n^2c\in N^2_0(c)$  be arbitrary.
Let $nc\in N(c)$ and  $c_i(nc)=(nc_x,nc_y).$
Define $c_i(n^2c)$ similar to that of $n^2a\in N^2_0(a),$ as follows.
\begin{myenu} [label= Case  \arabic*.]
\item If $n^2c\in \cup_{j<i}P_j$ and for some  $u,n^2c\in S_u,$
 we define $c_i(n^2c)$ as follows.
$$ \left\{\begin{array}{ll}
(-r(c),-r+2\delta (i+1))
& \mbox{if $nc_y\not =-r+2\delta (i+1),$  }\\
(-r/2+\delta(i+1),-r+2\delta(i+1)+r(c))
& \mbox{if $nc_y =-r+2\delta (i+1).$} \end{array}\right.$$

\item If $n^2c\in \cup_{j<i}P_j$ and for all  $u,n^2c\notin S_u,$
 we define $c_i(n^2c)=(0,-r+2\delta(i+1)).$

\item If $n^2c\notin \cup_{j<i}P_j,$   we define $c_i(n^2c)$ as follows.
$$ \left\{\begin{array}{ll}
(0,-r(c)) & \mbox{; $n^2ca\notin E(G)\nd n^2cb\notin E(G),$ }\\
(0,-r+2\delta (i+1)) &  \mbox{; $n^2ca\in E(G)\nd n^2cb\in E(G),$ }\\
(0,-r(c)+\delta ) & \mbox{; $n^2ca\in E(G)\nd n^2cb\notin E(G),$ }\\
(-r+2\delta(i+1)+r(c),-r+2\delta(i+1))   & \mbox{; $n^2ca\notin E(G)\nd n^2cb\in E(G).$} \end{array}\right.$$
\end{myenu}

\item Finally, let  $v\in V(G)\sm  V_i.$
We define $c_i(v)$ as follows.

If $v\in \cup_{j<i}P_j,$ define
$$c_i(v):=(-r(v)/2,-r(v)/2).$$
If $v\notin \cup_{j<i}P_j,$  we define $c_i(v)$ as follows.
$$ \left\{\begin{array}{ll}
(-r(v),-r(v)) & \mbox{; $va\notin E(G), vb\notin E(G), vc\notin E(G),  $ }\\
(-r(v),-r(v)+\delta)       & \mbox{; $va\in E(G), vb\notin E(G), vc\notin E(G),  $ }\\
(-r+2\delta(i+1)+r(v),-r+2\delta(i+1)+r(v))    & \mbox{; $va\notin E(G), vb\in E(G), vc\notin E(G),  $ }\\
(-r(v)+\delta,-r(v))       & \mbox{; $va\notin E(G), vb\notin E(G), vc\in E(G),  $ }\\
(-r(v)+2\delta,-r+2\delta(i+1)+r(v))      & \mbox{; $va\in E(G), vb\in E(G), vc\notin E(G),  $ }\\
(-r(v)+\delta,-r(v)+\delta)     & \mbox{; $va\in E(G), vb\notin E(G), vc\in E(G),  $ }\\
 (-r+2\delta(i+1)+r(v),-r(v)+2\delta)    & \mbox{; $va\notin E(G), vb\in E(G), vc\in E(G),  $ }\\
(-r(v)+2\delta,-r(v)+2\delta)     & \mbox{; $va\in E(G), vb\in E(G), vc\in E(G).$}
 \end{array}\right.$$
\end{myenu}

\subsection{Some Properties of the Super Triplet's  Case}

Before we move further, let us collect a few properties of the super triplet's case.

Recall the numbered cases and subcases in beginning of this section, while introducing the super triplet's case.
In the proof of each of the following propositions, we will use the same numbering in order to refer to the corresponding cases and subcases.

\begin{prop}\label{ssccond1}
   If $nb\notin \cup_{j<i}P_j, nba\notin E(G)\nd nbc\notin E(G),$ then \textit{$nb$ is not the  central vertex} of some star $S_u\cup \{u\}.$
\end{prop}
\begin{proof}
Assume that $nb$ is a central vertex of such a star. Then $b\in S_{nb}.$  The two cases of super triplets leads to contradictions, as can be seen below.
\begin{myenu}[label=Case \Roman*.]
\item
Since we picked $a,b,c,$ therefore $ac\notin E(G).$ Also $nba\notin E(G)\nd nbc\notin E(G).$ Therefore $a,nb,c$ form an independent set. As stated earlier $b\in S_{nb}$ therefore $a,nb,c$ should have been immediately picked in the previous step of picking $a,b,c.$ This is a contradiction.
\item
Since $nb\notin \cup_{j<i}P_j,$ by our Algorithm \ref{alg41},  $S_{nb}\subset \cup_{j<i}P_j.$
 But $b \in S_{nb}\nd b\in P_i.$ This is a contradiction.
\end{myenu}
Hence $nb$ is not a central vertex  of some star $S_u\cup \{u\}.$
\end{proof}

\begin{prop}\label{ssccond2}
Let $na\in N(a), c_i(na)=(na_x, na_y)\nd n^2a\in N^2_0(a).$ If
 $n^2a\in S_{na}, n^2a\in \cup_{j<i}P_j\nd na_x\ne -r+2\delta (i+1),$ then $n^2ab\notin E(G).$
\end{prop}
\begin{proof}
Consider the two cases of super triplets, as under.
\begin{myenu}[label=Case \Roman*.]
\item Consider the two subcases as follows.
\begin{myenu}[label=Subcase \roman*.]
\item In this case $b\in S_{nb},$ also $n^2a\in S_{na}.$ Therefore $n^2ab\notin E(G).$
\item In this case such $n^2a$ does not exist as $N^2(a)\sm \{a\}=\emptyset.$
\end{myenu}
\item As $S_{na}$ exists and $a\in P_i$ therefore $na\in \cup_{j<i}P_j$ (due to Algorithm \ref{alg41}.
Therefore $na_x =-r+2\delta (i+1),$ due to our assignment of dimensions in subsection \ref{assdimsec}.
\end{myenu}
Hence the result.
\end{proof}

\begin{prop}\label{ssccond3}
Let $nc\in N(c), c_i(nc)=(nc_x, nc_y)\nd n^2c\in N^2_0(c).$ If
$n^2c\in S_{nc}, n^2c\in \cup_{j<i}P_j\nd nc_y\ne -r+2\delta (i+1),$ then $n^2cb\notin E(G).$
\end{prop}
\begin{proof}
  Analogous to Proposition \ref{ssccond2}.
\end{proof}

  \begin{prop}\label{ssccond4}
  Let $P_k=\{a,b,c\},  u\in P_i\nd nu\in P_j\cap N(u)$
and let $u\in V(G)\sm  N_k.$
 If $k<min\{i,j\},$ then
  $m(u)>k.$
\end{prop}

\begin{proof}
Consider the two cases of super triplets, as under.
\begin{myenu}[label=Case \Roman*.]
\item
Assume that $m(u)\leq k.$
Then $m(u)\ne k,$ as $u\in V(G)\sm  N_k.$
Therefore  $m(u)<k.$
Also,  there exists some $w\in P_{m(u)}$
such that either $w\in N(u)\mbox{ or }w\in N^2(u).$

\hspace{.2in} If   $w\in N(u),$ then by construction (as $u\in P_i;i>k\nd nu\in P_j;j>k$), $u$ is the central vertex of some star and  $\{w,nu\}\subset S_u.$ Therefore there is some $w_1\in S_u$ such that
 $w_1\in P_k.$ Hence $u\in N(w_1)\nd  w_1\in P_k,$
a contradiction to the hypothesis that $u\in V(G)\sm  N_k.$ 

\hspace{.2in} If $w\in N^2(u),$ we write $w=nv,$ for some $v\in N(u).$
Then $v$ is the center of some star and $\{w,u\}\subset S_v.$
Therefore, there is some $w_2\in S_v$ such that $w_2\in P_k.$
  Hence $u\in N(v)\subset N^2(w_2)\nd  w_2\in P_k,$
a contradiction to the hypothesis that $u\in V(G)\sm  N_k.$ 
\item  There are three cases as under.
\begin{myenu}
\item[(a)]    \textbf{$u$ is a vertex of star,}   say $u\in \{w\}\cup S_w.$ Assume that $m(u)\leq k.$
Note that in this case, either $u$ is picked in some $P_i$ with $i\leq k,$ or the whole of $S_w$ is picked before $P_k.$
Therefore, either $i\leq k\mbox{ or }j\leq k,$ a contradiction.
%
\item[(b)] \textbf{$u$ is a vertex of matching.} In this case, clearly $m(u)>k$ as there are only two vertices $v\nd nu$ under consideration.
\item[(c)] \textbf{$u$ is a vertex of triangle.} In this case, the pseudo-neighbors of the triangle are defined in such a way that $m(u)>k.$
\end{myenu}
\end{myenu}
This establishes the result.
\end{proof}
\subsection{$G$ is Realizable as a SIG via Our Mapping}
By now, we have mapped the vertex set corresponding to the super triplet's case, on an Euclidean space by assigning the coordinates with respect to each set of vertices picked in the algorithm \ref{alg41}.
For convenience, we will use the same symbol $v$ for the image of $v,$ under this mapping.

In this sense, the vertex set $V(G)$ is now projected in an Euclidean space endowed with sup metric.
We now prove that the SIG of this mapped vertex set is isomorphic to the given graph.
We prove our main result through a series of lemmas.

In these lemmas, for $u, v\in G$,  we will use the notation $|c_k(u)-c_k(v)|,$ even when $c_k$ represents a pair of Euclidean dimensions. In that case, as an abuse of notation,  it will represent the sup-norm in those two dimensions.

\subsection{Proof of (\ref{LemEq1})}
\begin{lemma}\label{ssclem1}
If $P_k$ is a super triplet, $u\in V(G) \nd nu\in N(u),$ then
\begin{equation}\label{ssclem1eqmaain}
|c_k(u)-c_k(nu)| \leq r(u).
\end{equation}
That is, $u\nd nu$ are at a distance of $r(u)$ or less in the plane.
\end{lemma}

\begin{proof}
Write $P_k=\{a,b,c\},$ where $a,b,c$ are as discussed in the introduction of super triplets, in the beginning of this section.
Note that there exist $i\nd j$ such that $u \in P_i\nd nu \in P_j.$
Then, by the definitions of $m(u)\nd r(u),$ we have
$$m(u)\leq \min\{i,j\}\nd  r(u)=r-2\delta m(u) \geq r-2\delta \min\{i,j\}.$$


Note that for each $u\in V(G),$ the set $N(u)\cup N^2(u)$ is closed under the pseudo-neighborhood operation.
Thus,  it is enough to show that
$|c_k(nu)-c_k(n^2u)|\leq r(u),$
for each $nu\in N(u)$ and  $n^2u\in N^2(u).$ Hence we consider the following cases.
\begin{myenu}[label=Case \arabic*.]
\item\label{ssclem1case1} $\bm{u\in N(a)}.$ Then $m(u)\leq k.$
Consider the following sets.
\begin{align*}
&A_1:= \big{\{} (-r(u),0),(-r(u)+\delta ,0) \big{\}}, \\
&A_2:=\big{\{} (-r+2\delta (k+1),0), (-r+2\delta(k+1),-r+2\delta(k+1)+r(u))\big{\}}, \\
&A_3:=\big{\{}(-r+2\delta(k+1)-r(u), r(u))\big{\}},\\
&A_4:=\big{\{}(-r+2\delta (k+1), -r(u))\big{\}},\\
&A_5:=\big{\{}(-r+2\delta(k+1)+r(u),-r/2+\delta(k+1))\big{\}}.
\end{align*}
Since $u\in N(a),$ observe that either $c_k(u)\in A_1\mbox{ or }c_k(u)\in A_2.$ \\
If $c_k(u)\in A_1,$ then $c_k(nu) \in \cup_{p=1}^4 A_p.$ If $c_k(u)\in A_2,$ then $c_k(nu)\in \cup_{p\ne 4}A_p.$
Hence,  we observe that (\ref{ssclem1eqmaain}) holds if each of
\begin{gather*}
|2\delta (k+1)-r|, |2\delta (k+1)-r+r(u)|, |-r/2+\delta (k+1)|,\\
|\delta (2k+1)-r|, |\delta (2k+1)-r+r(u)|
\end{gather*}
is bounded by $r(u).$ This holds, due to the fact that $m(u)\leq k$ and by the definition of $\del.$
Hence inequality  (\ref{ssclem1eqmaain}) is satisfied for this case.
\item \bm{$u\in N(b).$} $m(u)\leq k.$ Write
\begin{align*}
&\alpha_1:=(\delta, \delta),\\
&\alpha_2:=(\delta, -r+2\delta(k+1)+r(u)),\\
&\alpha_3:=(-r+2\delta(k+1)+r(u), \delta),\\
&\alpha_4:=(\delta+r(u), \delta+r(u)),\\
&\alpha_5:=(\delta -r(u), \delta -r(u)),\\
&\alpha_6:=(\delta -r(u), -r+2\delta(k+1)),\\
&\alpha_7:=(-r+2\delta(k+1), \delta -r(u)),\\
&\alpha_8:=(-r+2\delta(k+1)+r(u),-r+2\delta(k+1)+r(u)).
\end{align*}
Now, we observe the following.\\
If $c_k(u)=\alpha_1,$ then $c_k(nu)\in \{\alpha_1,\alpha_2,\alpha_3,\alpha_4,\alpha_5,\alpha_8\}.$\\
 If  $c_k(u)=\alpha_2,$ then $c_k(nu)\in \{\alpha_1,\alpha_2,\alpha_3,\alpha_4,\alpha_6,\alpha_8\}.$\\
  If  $c_k(u)=\alpha_3,$ then $c_k(nu)\in \{\alpha_1,\alpha_2,\alpha_3,\alpha_4,\alpha_7,\alpha_8\}.$ \\
If $c_k(u)=\alpha_8,$ then $c_k(nu)=\alpha_4,$  due to Proposition \ref{ssccond1}.\\
In each of these cases, as above, one can verify  that (\ref{ssclem1eqmaain}) holds.
\item \bm{$u\in N(c).$} This case is analogous to \ref{ssclem1case1}.

\item \label{case 1.1.4.-ssc}
 $\bm{u\in V(G)\sm  N_k.}$ Then  $nu\in V(G)\sm  N_k.$
\begin{myenu}
\item \label{1.4.1-ssc} \bm{$k< min\{i,j\}.$} By Proposition \ref{ssccond4},  we have $m(u) > k.$
%
Note that
each  of $c_k(u)[0],$ $c_k(nu)[0],$ $c_k(u)[1]\nd c_k(nu)[1]$ belongs to the set
$$\big{\{}-r(u),-r(u)+\delta ,-r(u)+2\delta ,-r+2\delta (k+1)+r(u)\big{\}}.$$

Therefore, inequality (\ref{ssclem1eqmaain}) holds,
if and only if we have
\begin{gather*}
2\del<r(u),  |2r(u)-r+2\delta (k+1)|\leq r(u), \\
|2r(u)-r+ \delta (2k+1)|\leq r(u)
\mbox{ and }|2r(u)-r+2\delta k)|\leq r(u).
\end{gather*}
As earlier, we only need to show that
\begin{align*}
r(u)&\leq r-2\delta (k+1) \leq 3r(u),\\
r(u)&\leq r-\delta (2k+1) \leq 3r(u)\\
\mbox{ and } r(u)&\leq r-2\delta k \leq 3r(u).
\end{align*}
The first parts of these three inequalities hold due to $m(u) > k.$
The second parts of these inequalities hold due to our choice of $\del.$
Hence, for this case, inequality (\ref{ssclem1eqmaain}) is verified.
\item \bm{$k > min\{i,j\}.$}
In this case, note that $c_k(u)[0]=c_k(u)[1]=- r(u)/2$ and both $c_k(nu)[0]$ and $c_k(nu)[1]$ belong to set
$$  \big{\{}-\frac{r(u)}{2},-r(u),-r(u)+\delta, -r(u)+2\delta ,-r+2\delta (k+1)+r(u)\big{\}}.$$
Further, note that inequality (\ref{ssclem1eqmaain}) holds if and only if each of
$|\delta-r(u)/2|,  |2\delta-r(u)/2|, |2\delta(k+1)-r+ 3r(u)/2|$ is bounded above by $r(u).$
Since $-r+ 2r(u)=r-4\delta m(u) =r(1-\frac{2}{3n}m(u))\geq 0, $ \wo
\begin{equation}\label{nfjsafhfdssadas}
 \delta-\frac{r(u)}{2} \leq 2\delta-\frac{r(u)}{2} \leq 2\delta(k+1)-r+ \frac{3}{2}r(u).
 \end{equation}
\indent Now, \nt $2\delta(k+1)-r+ 3r(u)/2\leq r(u)$ if
$2\delta(n+1)+ r(u)/2\leq r,$  that is,  $4\delta(n+1)+ r(u)\leq 2r$  or
$4\delta(n+1)-2\delta m(u)\leq r.$  This holds if
$4\delta(n+1)\leq r,$ which is true due to our choice of $\del.$

\indent Similarly, the inequality $\delta-r(u)/2\geq -r(u)$ holds, again due to our definition of $\del.$
Along with (\ref{nfjsafhfdssadas}), these imply  (\ref{ssclem1eqmaain}) for this case.
%
%
%

\end{myenu}
\end{myenu}
Therefore in all cases, the inequality (\ref{ssclem1eqmaain}) is verified. Hence the result.
\end{proof}

\subsection{Proof of (\ref{LemEq2})}

Similarly,   verifying step by step, we can establish the following.
\begin{lemma}\label{ssclem2}
Let $\{u,v\}\subset V(G)$  such that $u\ne v, u \in P_i, v \in P_j, i\leq j\nd P_j$ is a super triplet. Then
\begin{equation}\label{ssclem2eqmaain}
|c_j(u)-c_j(v)| \geq \max \{r(u), r(v)\}.
\end{equation}
\end{lemma}
 \begin{proof}
By hypothesis, \wv $m(u)\leq i$ and therefore $r(u)=r-2\del m(u) \geq r-2\del i.$
 Write $P_j:=\{a,b,c\}.$  Before we establish  (\ref{ssclem2eqmaain}),  case by case, we present a few inequalities which will be used in this proof. These can be established from the definitions of $\del, r(u)\nd r(v).$
 \begin{align*}
|    \delta +r(v)+r-2\delta (j+1)+r(u)|&=3r-\del(1+6n)>r,\\
|    \delta +r(v)+r-2\delta (j+1)|&=2r-\del(1+2j+2m(v))>r,\\
|    \delta +r(v)-\del+ r(u)|&=2r-2\del(m(u)-m(v))>r,\\
|    \delta +r(v)+r(u)|&=2r+\del(1-2m(v)-2m(u))>r,\\
|    \delta +r(v) + r(u)/2|&=3r/2+\del(1-2m(v)-m(u  ))>r,\\
\nd |    \delta +r(v)+(r-2\delta (j+1))/2|&=3r/2-\del(j+2m(v))>r.
 \end{align*}
Now, we consider the following cases.
\begin{myenu}[label=Case \arabic*.]
\item\label{ssclem2case1}
\bm{$v=a.$} Then $c_j(v)=(-r+2\delta (j+1)-r(v),r(v))$ and $c_j(u)$ belongs to the set
\begin{align*}
  \big\{
(\delta&+r(u), \delta+r(u)),(r(u),-r+2\delta(j+1)-r(u)),(-r+2\delta(j+1),0),(\delta,\delta),\\
&(0,-r+2\delta(j+1)), (-r+2\delta(j+1),-r(u)),(-r+2\delta(j+1)+r(u),-r/2+\delta(j+1)),\\
&(\delta-r(u), \delta-r(u)), (\delta-r(u), -r+2\delta(j+1)),(-r+2\delta(j+1), \delta-r(u)),\\
&(-r(u),-r+2\delta(j+1)),(-r/2+\delta(j+1), -r+2\delta(j+1)+r(u)),(-r(u)/2,-r(u)/2)\big\}.
\end{align*}
\\Hence (\ref{ssclem2eqmaain}) holds true.
\item
\bm{$v=b.$} Then {$c_j(v)= (\delta +r(v),\delta +r(v))$} and $c_j(u)$ belongs to the set
 \begin{align*}
    \big\{
(-r+&2\delta(j+1)-r(u),r(u)),(r(u),-r+2\delta(j+1)-r(u)),(-r+2\delta(j+1),0),(\delta,\delta),\\
&(0,-r+2\delta(j+1)),(-r+2\delta(j+1),-r(u)),(-r+2\delta(j+1)+r(u),-r/2+\delta(j+1)),\\
&(\delta-r(u), \delta-r(u)),(\delta-r(u), -r+2\delta(j+1)),(-r+2\delta(j+1), \delta-r(u)),\\
&(-r(u),-r+2\delta(j+1)),(-r/2+\delta(j+1), -r+2\delta(j+1)+r(u)),(-r(u)/2,-r(u)/2)\big\}.
 \end{align*}
\\Hence (\ref{ssclem2eqmaain}) holds true.
\item \bm{$v=c.$}
This case is analogous to \ref{ssclem2case1}
\end{myenu}
This establishes the lemma.
\end{proof}

This lemma, along with the previous modules, implies that
$$r_v = r(v)\foa v\in V(G).$$

\subsection{Proof of (\ref{LemEq3})}

To conclude that our mapping is a sphere of influence graph, we need to establish  a few more lemmas as follows.

\begin{lemma}\label{ssclem3}
Let $\{u, v\}\subset V(G)$ \st $ u\ne v, uv\notin E(G)\nd  \{u,v\}\subseteq S_w,$ for some $w\in V(G).$
If $u \in P_i, v \in P_j$ with $i\leq j $ and
 $P_j$ is a perfect triplet, then
\begin{equation}\label{ssclem3maaineq}
|c_j(u)-c_j(v)| \geq r(u)+r(v).
\end{equation}
\end{lemma}
\begin{proof}
\Nt by hypothesis, we have $$r(u) \geq r-2\delta i \mbox{ and } r(v) \geq r-2\delta j.$$
Write $P_j=\{a,b,c\}$ and consider the following cases.
\begin{myenu}[label=Case \arabic*.]
\item\label{ssclem3case1} \bm{$v=a.$} Then $c_j(v)=(-r+2\delta (j+1)-r(v),r(v))\nd c_j(u)$ belongs to
$$\big\{(-r+2\delta (j+1),-r(u)),(-r+2\delta (j+1)+r(u),-r/2+\delta (j+1))\big\}.$$
\item \bm{$v=b.$} Then $c_j(v)=(\delta +r(v),\delta +r(v))\nd c_j(u)$ belongs to the set
$$ \big\{(\delta -r(u),\delta -r(u)),(\delta -r(u),-r+2\delta (j+1)),(-r+2\delta (j+1),\delta -r(u))\big\}.$$
\item \bm{$v=c.$} Similar to \ref{ssclem3case1}\qedhere
\end{myenu}
\end{proof}
\subsection{Proof of (\ref{LemEq4})}

\begin{lemma}\label{ssclem4}
Let $\{u, v\}\subset V(G)$ \st $ u\ne v, uv\notin E(G)\nd  \{u,v\}\not \subseteq S_w,$ for any $w\in V(G).$
If $u \in P_i, v \in P_j$ with $i\leq j $, and
 $P_i$ is a perfect triplet, then
\begin{equation}\label{ssclem4maaineq}
|c_i(u)-c_i(v)| \geq r(u)+r(v).
\end{equation}
\end{lemma}

\begin{proof}
\Nt by hypothesis, we have $$r(u) \geq r-2\delta i \mbox{ and } r(v) \geq r-2\delta j$$
Write $P_i=\{a,b,c\}$ and consider the following cases.
\begin{myenu}[label=Case \arabic*.]
\item\label{ssclem4case1}  \bm{$u=a.$} Then $c_i(u)=(-r+2\delta(i+1)-r(u), r(u))$ and   $c_i(v)$ belongs to the set
\begin{align*}
 &\big{\{}(\delta+r(v), \delta+r(v)),(r(v),-r+2\delta(i+1)-r(v)), (-r+2\delta(i+1)+r(v), \delta),\\
&(0,-r(v)),(-r+2\delta(i+1)+r(v),-r+2\delta(i+1)),(-r+2\delta(i+1)+r(v),\\
&-r+2\delta(i+1)+r(v)),(-r(v),-r(v)),(-r(v)+\delta,-r(v)),\\
&(-r+2\delta(i+1)+r(v),-r(v)+2\delta)\big{\}}.
\end{align*}
\item

    \bm{$u=b.$} Then $c_i(u)=(\delta+r(u), \delta+r(u))$ and both $c_i(v)$ belongs to set
\begin{align*}
\big{\{}
(-r+2&\delta(i+1)-r(v),r(v)),(r(v),-r+2\delta(i+1)-r(v)),(-r(v),0), (-r(v)+\delta,0), \\
&(0,-r(v)), (0,-r(v)+\delta), (-r(v),0),(-r(v)+\delta, 0), (-r(v),-r(v)), \\
&(-r(v),-r(v)+\delta),  (-r(v)+\delta,-r(v)), (-r(v)+\delta,-r(v)+\delta)\big{\}}.
\end{align*}
    \item \bm{$u=c.$}  Similar to \ref{ssclem4case1}\qedhere
 \end{myenu}
\end{proof}

\subsection{Proof of (\ref{LemEq5})}
\begin{lemma}\label{ssclem5}
If  $P_k$ is a super triplet and   $\{u, v\}\subseteq V(G)$ such that $uv\in E(G),$ then
  \begin{equation}\label{ssclem5maaineq}
 |c_k(u)-c_k(v)| < r(u)+r(v).
\end{equation}
\end{lemma}

\begin{proof}
Write $P_k=\{a,b,c\}.$
Suppose that $u \in P_i\nd v \in P_j.$ Then
  $$r(u) \geq r-2\delta i \mbox{ and } r(v) \geq r-2\delta j.$$
We will establish (\ref{ssclem5maaineq}) case by case as follows.
\begin{myenu}[label=Case \arabic*.]
\item \bm{$k\notin \{i,j\}.$} Let $Q:=(-\delta ,-\delta )$ and
\begin{align*}
A(u) = & \big{\{}-r+2\delta(k+1),-r(u),-r(u)+\delta ,\delta,-r+2\delta(k+1)+r(u),0,\\
&\delta -r(u),-r/2+\delta(k+1),-r(u)/2,-r(u)+2\delta,\big{\}}.
\end{align*}
Both $c_k(u)[0]\nd c_k(u)[1]$ belong to $A(u)$ and both $c_k(v)[0]\nd c_k(v)[1]$ belong to $A(v).$
$$|c_k(u)-Q|<r(u)\nd |c_k(v)-Q|<r(v).$$
As usual, the triangle inequality implies the required inequality.
\item \bm{$k\in \{i,j\}.$} Without loss of generality, we assume that $k=i.$
\begin{myenu}
\item\label{ssclem5case1point1} \bm{$u=a.$} Then
$c_k(u)=(-r+2\delta(k+1)-r(u), r(u))$ and
\begin{align*}
c_k(v)[0]\in \big{\{}-r+&2\delta(k+1),-r(v),-r(v)+\delta ,\delta,0,\delta -r(v),\\
&-r/2+\delta(k+1),-r(v)/2,-r(v)+2\delta\big{\}}.
\end{align*}
\begin{align*}
c_k(v)[1]\in \big{\{}0,-r+&2\delta(k+1)+r(v), \delta,-r+2\delta(k+1),-r(v)+\delta,\\
& \delta -r(v),-r(v)/2,-r(v)+2\delta)\big{\}}.
\end{align*}

Let $Q:=(-r+2\delta (k+1)-\frac{\delta }{2},\frac{\delta }{2}).$
$$|c_k(u)-Q|<r(u)\nd |c_k(v)-Q|<r(v).$$
Hence, by the triangle inequality, \wo  the required inequality.
\item \bm{$u=b.$} then
$c_k(u)=(\delta+r(u), \delta+r(u))$ and
$$A=\bigg{\{}-r+2\delta(k+1),\delta,-r+2\delta(k+1)+r(v),0,-r(v)/2,-r(v)+2\delta,\frac{-r+2\delta(k+1)}{2}\bigg{\}}.$$
Applying Proposition \ref{ssccond2}, $c_k(v)[1]\neq -r(v)$ and its analogue Proposition \ref{ssccond3}, we obtain that
$c_k(v)[0]\neq -r(v).$
Therefore, both $c_k(v)[0]\nd c_k(v)[1]$ belong to the set $A.$
Let $Q:=(\frac{3\delta }{2},\frac{3\delta }{2}).$
$$|c_k(u)-Q|<r(u)\nd |c_k(v)-Q|<r(v).$$
As usual, the triangle inequality implies the required inequality.
\item \bm{$u=c.$} Similar to \ref{ssclem5case1point1}\qedhere
\end{myenu}
 \end{myenu}
\end{proof}


%

\section{\bf Proof of The Main Result}

\noindent Using all the previous modules, we obtain the following result.
\begin{theorem}
For $v_1, v_2\in V(G),$
$$v_1v_2\in E(G)\mbox{  if and only if }\rho(v_1, v_2) < r_{v_1}+r_{v_2}.$$
\end{theorem}
Therefore the sphere of influence graph of our mapping of $V(G)$ on an Euclidean space is isomorphic to $G.$
In other words, $G$ is realizable in a Euclidean space, whose dimension is fixed according to our algorithm.
Next, we   estimate a bound on the dimension of this Euclidean space.\\

\noindent{\bf Proof of Theorem \ref{MainThmSIG3}.}
Recall that our $G$ is any graph of order $n$ and has no  isolated vertex. We need to establish that
\begin{equation}\label{MainThmMainEq}
SIG(G)\leq \bigg{ \lfloor}\frac{2n}{3}\bigg{ \rfloor}+2.
\end{equation}
 Let  $\mathcal{C}_t\nd \mathcal{C}_p,$ respectively, denote the collections of  triplets and  pairs of vertices of $G,$
 picked up as per Algorithm \ref{alg41}, which do not fall in the category of residual sets or random pairs of module I and II.
 This algorithm also ensures that out all the three possibilities of  our residual sets,
 at most one is picked up. If no residual set  is picked up, define $R_t=\emptyset.$ Otherwise, let $R_t$ denote the residual set, which has been picked up.

Further, note that there are four possibilities for the random pairs or singleton sets to be picked up, which occur at
steps \ref{CharVerS18}-\ref{CharVerS19}, \ref{CharVerS30}, \ref{CharVerS37}  or \ref{CharVerS45} of   Algorithm \ref{alg41}.
However, both of the  steps \ref{CharVerS37}  and \ref{CharVerS45}  can't occur together in  Algorithm \ref{alg41}.
Therefore, there  are at most three possibilities for the random pairs or singleton sets to be picked up.

%
Thus,   there exist $k_1,k_2, k_3\in \{0,1,2\}$ such that
\begin{equation}\label{MainThmEq1}
  n=3|\mathcal{C}_t|+2|\mathcal{C}_p|+|R_t| +k_1+k_2+k_3.
  \end{equation}

Recall  our embedding of $G$ into an Euclidean space, as a sphere of influence graph. We assigned
$\lceil\log_2(m+1)\rceil $ dimensions for the residual set of size $m$ and $k_1, k_2, k_3$ extra dimensions for each of random pairs and singleton sets of sizes $k_1, k_2\nd k_3.$ For other pairs and triplets, we assigned
two dimensions per triplet and one dimension per pair.  Therefore,  there exist $k_1,k_2, k_3\in \{0,1,2\}$ such that
\begin{equation}\label{MainThmEq2}
SIG(G) \leq 2|\mathcal{C}_t|+|\mathcal{C}_p|+\lceil\log_2(|R_t|+1)\rceil+k_1+k_2+k_3.
  \end{equation}

Now, we divide the proof into two cases. First we consider the case when a set of vertices has been picked up using step \ref{CharVerS30} of   Algorithm \ref{alg41}. In this case, no residual set is picked up and (\ref{MainThmEq1}) and (\ref{MainThmEq2}) become
\begin{equation}\label{MainThmEq1A}
  n=3|\mathcal{C}_t|+2|\mathcal{C}_p|+k_1+k_2+k_3.
  \end{equation}
\begin{equation}\label{MainThmEq2A}
\nd SIG(G) \leq 2|\mathcal{C}_t|+|\mathcal{C}_p|+ k_1+k_2+k_3.
  \end{equation}
Substituting  (\ref{MainThmEq1A})  in  (\ref{MainThmEq2A}), we obtain
\begin{align*}
SIG(G)
&\leq  2|\mathcal{C}_t|+|\mathcal{C}_p|+ k_1+k_2+k_3\\
&\leq \bigg\lfloor \frac{2n-(k_1+k_2+k_3)}{3} \bigg\rfloor   +k_1+k_2+k_3 \\
&=\bigg\lfloor \frac{2n-(k_1+k_2+k_3)}{3} + k_1+k_2+k_3  \bigg\rfloor   \\
&=\bigg\lfloor \frac{2n}{3}+\frac{k_1+k_2+k_3}{3} \bigg\rfloor   \\
&\leq \bigg\lfloor \frac{2n}{3}+ 2\bigg\rfloor   \\
&\leq \bigg{ \lfloor}\frac{2n}{3}\bigg{ \rfloor}+2.
\end{align*}

Now, consider the   case when no set of vertices has been picked up using step \ref{CharVerS30} of   Algorithm \ref{alg41}. In that case,  there are at most two random pairs or singletons and at most one residual set is picked up.
In this case (\ref{MainThmEq1}) and (\ref{MainThmEq2}) become
\begin{equation}\label{MainThmEq1B}
  n=3|\mathcal{C}_t|+2|\mathcal{C}_p|+|R_t|+k_1+k_2.
  \end{equation}
\begin{equation}\label{MainThmEq2B}
\nd SIG(G) \leq 2|\mathcal{C}_t|+|\mathcal{C}_p|+\lceil\log_2(|R_t|+1)\rceil+k_1+k_2.
  \end{equation}
   Routine calculus and number theoretic arguments can be applied to prove  that every non negative integer $k$ satisfies the relations
\begin{equation}\label{MainThmEqCal}
\bigg\lceil \frac{2k}{3}\bigg\rceil=\bigg\lfloor \frac{2(k+1)}{3}\bigg\rfloor
\nd  \lceil\log_2(k+1)\rceil\leq   \bigg\lceil \frac{2k}{3} \bigg\rceil.
    \end{equation}
Applying (\ref{MainThmEq1B}) and (\ref{MainThmEqCal})  in  (\ref{MainThmEq2B}), we obtain
\begin{align*}
SIG(G)
&\leq  2|\mathcal{C}_t|+\frac{4|\mathcal{C}_p|}{3}+\bigg\lceil\frac{2|R_t|}{3}\bigg\rceil +k_1+k_2 \\
&\leq  \bigg\lceil 2|\mathcal{C}_t|+\frac{4|\mathcal{C}_p|}{3}+\frac{2|R_t|}{3}\bigg\rceil+k_1+k_2 \\
&= \bigg\lceil \frac{2(n-k_1-k_2)}{3}\bigg\rceil+k_1+k_2 \\
&= \bigg\lfloor \frac{2(n-k_1-k_2+1)}{3}\bigg\rfloor+k_1+k_2 \\
&= \bigg\lfloor \frac{2(n-k_1-k_2+1)}{3}+k_1+k_2\bigg\rfloor  \\
&= \bigg\lfloor \frac{2n}{3}+\frac{k_1+k_2+2}{3}\bigg\rfloor \\
&\leq \bigg\lfloor \frac{2n}{3}+2\bigg\rfloor =\bigg{ \lfloor}\frac{2n}{3}\bigg{ \rfloor}+2.
\end{align*}
This proves (\ref{MainThmMainEq}) in both cases.
Further, if $n$ is not a multiple of $3,$ then
$$SIG(G)\leq \bigg{ \lfloor}\frac{2n}{3}\bigg{ \rfloor}+2=\bigg{ \lceil}\frac{2n}{3}\bigg{ \rceil}+1.$$
The result is established. \null\hfill\qedsymbol

Whether the conjecture is true or not,  remains an open problem.

\end{document}